[1994/12/01]
\documentclass[10pt, reqno]{amsart}
\usepackage[english, activeacute]{babel}
\usepackage{amsmath,amsthm,amsxtra}
\usepackage{epsfig}
\usepackage{amssymb}
\usepackage{latexsym}
\usepackage{amsfonts}
\usepackage{hyperref}
\usepackage{pdftricks}
\usepackage{wrapfig}
\usepackage{tikz}
\usetikzlibrary{arrows,calc,decorations.pathreplacing}
\definecolor{light-gray1}{gray}{0.90}
\definecolor{light-gray2}{gray}{0.80}
\definecolor{light-gray3}{gray}{0.60}

\title{Stability of integrable and nonintegrable structures}
\author{Claudio Mu\~noz}
\address{CNRS and Laboratoire de Math\'ematiques d'{}Orsay UMR 8628, B\^at. 425 Facult\'e des Sciences d'{}Orsay,
Universit\'e Paris-Sud F-91405 Orsay Cedex France}
\email{claudio.munoz@math.u-psud.fr}
\date{\today}
\subjclass[2000]{Primary 35Q51, 35Q53; Secondary 37K10, 37K40}
\keywords{generalized KdV equations, solitons, multi-solitons, stability, integrability, collision, breathers}
\thanks{This paper is part of a course taught for a graduate class at IMPA (Brazil). I would like to thank to Felipe Linares and IMPA's members for the kind invitation and its hospitality.}

\chardef\bslash=`\\ 

\newtheorem{thm}{Theorem}[section]
\newtheorem{cor}[thm]{Corollary}
\newtheorem{lem}[thm]{Lemma}
\newtheorem{prop}[thm]{Proposition}

\newtheorem{defn}[thm]{Definition}

\theoremstyle{remark}
\newtheorem{rem}{Remark}[section]

\newtheorem{Cl}{Claim}

\numberwithin{equation}{section}

\newcommand{\R}{\mathbb{R}}

\newcommand{\Z}{\mathbb{Z}}

\newcommand{\Com}{\mathbb{C}}
\newcommand{\la}{\lambda}

\newcommand{\al}{\alpha}
\newcommand{\bt}{\beta}
\newcommand{\ga}{\gamma}

\def\bm{\left( \begin{array}{cc}}
\def\endm{\end{array}\right)}

 \providecommand{\abs}[1]{\lvert#1 \rvert}

\newcommand{\ve}{\varepsilon}

\newcommand{\be}{\begin{equation}}
\newcommand{\ee}{\end{equation}}
\newcommand{\ba}{\left(\begin{array}{c}}
\newcommand{\ea}{\end{array}\right)}
\newcommand{\bea}{\begin{eqnarray}}
\newcommand{\eea}{\end{eqnarray}}
\newcommand{\bee}{\begin{eqnarray*}}
\newcommand{\eee}{\end{eqnarray*}}
\newcommand{\ben}{\begin{enumerate}}
\newcommand{\een}{\end{enumerate}}
\newcommand{\nonu}{\nonumber}

\newcommand{\eval}[2][\right]{\relax
  \ifx#1\right\relax \left.\fi#2#1\rvert}

\let\abs=\envert

\begin{document}
\begin{abstract}
In this paper we give a comprehensive account of several recent results on the stability of nontrivial soliton structures for some well-known non periodic dispersive models. We will focus on the simpler case of the generalized Korteweg-de Vries equations, covering the classical stability results by Bona, Souganidis, Strauss until the results by Martel and Merle and our recent collaborations with Miguel Alejo and Luis Vega.   
\end{abstract}
\maketitle \markboth{Stability of soliton structures} {Claudio Mu\~noz}
\renewcommand{\sectionmark}[1]{}
\tableofcontents

\newpage

\section{Introduction. Classical stability of solitons}

\medskip

\subsection{Introduction} 

Consider the generalized Korteweg-de Vries (gKdV) equation on the real line 
\be\label{gKdV}
u_t+(u_{xx} + u^p)_x =0,
\ee
($p$ is a positive integer) where $u=u(t,x)$ is a real-valued function, and $(t,x)\in \R^2$. The case $p=2$ is the famous  Korteweg-de Vries (KdV) equation, a model for shallow water waves and a canonical example of a dispersive nonlinear evolution equation. We assume $p$ integer since $u$ may have no definite sign, otherwise we must consider the odd power case $|u|^{p-1}u$, $p\geq 1.$

\medskip

Even if \eqref{gKdV} is a simple 1d model, its study proved to be a very difficult problem (see e.g. the monograph by Linares and Ponce \cite{LP}). In particular, if one considers the Cauchy problem
\be\label{gKdVCP}
u_t+(u_{xx} + u^p)_x =0, \quad u(t=0)=u_0,
\ee
Kato first, and then Kenig, Ponce and Vega \cite{KPV} showed global existence and well-posedness if $u_0\in H^1(\R)$, for $p=2,3$ and $4$. These results have been improved by many others, see e.g. the works by Bourgain \cite{Bo1} the I-team \cite{CKSTT}.   For the purposes of this course, we only need the $H^1$ global well-posedness, that is

\begin{thm}[\cite{KPV}] 

Assume $u_0\in H^1(\R)$, with $p=2,3$ or $4$ in \eqref{gKdVCP}. Then there exists a unique  $u\in C(\R,H^1(\R))$ solution of \eqref{gKdVCP} in the Duhamel sense:
\[
u(t) = S(t)u_0 - \int_0^t S(t-s)[(u^p)_x(s)]ds, \quad S(t) := e^{-t \partial_x^3}.
\]
Moreover, one has
\[
\sup_{t\in\R} \|u(t)\|_{H^1(\R)} \leq C(\|u_0\|_{H^1(\R)}).
\]
\end{thm}  

Whenever $p\geq 5$, problem \eqref{gKdVCP} becomes a very difficult one: in particular, if $p=5$ it is well-known that there are blowing-up solutions (see Martel and Merle \cite{MMblow}). If $p>5$, very little is known, mainly because the problem is $L^2$ ``supercritical''. We will discuss some of these affirmations below. 

\medskip

\subsection{Solitons} 

In addition to the previous global well-posedness result, \eqref{gKdVCP} is interesting by its \emph{soliton} solutions. A soliton for \eqref{gKdVCP} is a solution of the form
\be\label{S}
u(t,x) = Q_c(x-ct-x_0), \quad c>0, \, x_0\in \R.
\ee
Note that $Q_c$ is a fixed profile depending only on $c>0$ (usually referred as the scaling, or the velocity of the soliton), and $x_0$ is a free shift parameter. Replacing \eqref{S} in \eqref{gKdV}, and assuming that $Q_c$ vanishes at infinity, we find that $Q_c$ must satisfy the equation
\[
(Q_c''-cQ_c+Q_c^p)' =0, 
\]
or 
\be\label{S1}
Q_c''-cQ_c+Q_c^p=0 \quad \hbox{ in } \ \R.
\ee
The solutions to this elliptic equation (actually, it is an ODE), are well-understood: $Q_c=Q_c(s)$ is given by the explicit formula
\[
Q_c(s) =c^{1/(p-1)} Q(\sqrt{c} s),
\] 
where\footnote{Note that $Q_c$ is even.}
\[
Q(s) \ (=Q_{c=1}) := \left( \frac{p+1}{2\cosh^2(\frac{(p-1)}{2}s)}\right)^{1/(p-1)}.
\]

\medskip

In other words, a soliton is a traveling wave solution with positive speed $c>0$. Note that the bigger $c>0$ is, the faster the soliton travels.

\medskip

We would like to emphasize that solitons are purely nonlinear objects: the associated linear equation, denoted as the Airy equation
\[
u_t + u_{xxx} =0,
\] 
has no soliton solutions.

\medskip

A very nice relation between solitons and the global well-posedness result can be established using the $L^2$ norm of each soliton. We have
\be\label{crit}
\|Q_c\|_{L^2(\R)} \sim c^{\frac1{p-1}-\frac14},
\ee
so that $\frac1{p-1}-\frac14>0$ if $p=2,3,4$ ($L^2$-subcritical regime), $\frac1{p-1}-\frac14=0$ if $p=5$ ($L^2$-critical regime), and $\frac1{p-1}-\frac14<0$ if $p>5$ ($L^2$-supercritical regime). 

\medskip

In other words, in a supercritical problem, small solitons (in the sense of the $L^\infty$ norm), are very large in the $L^2$ norm. Whenever $p=5$, all solitons have the same size. We can summarize these properties by saying that 
\[
\partial_c \int_\R Q_c^2 \Big|_{c=1} = : \int_\R \Lambda Q Q \  \begin{cases} >0, & p=2,3,4, \\ =0, & p=5, \\ <0, & p>5. \end{cases}
\]
Here we have denoted by $\Lambda Q_c$ the so-called scaling direction, namely
\be\label{LaQ}
\Lambda Q_c (s):= \partial_c Q_c(s) =\frac1c\Big( \frac{1}{p-1}Q_c(s) +\frac 12 s Q_c'(s) \Big),
\ee
and $\Lambda Q := \Lambda Q_c \Big|_{c=1}.$ Note that $\Lambda Q_c$ is an even function\footnote{In fact, $\Lambda Q_c$ is a very important direction for the dynamics of nontrivial perturbations of a soliton solution, for instance in the critical case $p=5$.} which satisfies the equation (just take derivative with respect to $c$ in \eqref{S1})
\be\label{LaQc}
\mathcal L \Lambda Q_c = -Q_c.
\ee

\subsection{Stability of solitons}

In this paragraph we will discuss the stability problem for soliton solutions. Assume that $u_0\in H^s(\R)$ satisfies 
\be\label{1p0}
\|u_0-Q_c\|_{H^s} <\alpha,
\ee
where $\alpha\ll 1$ and $s\geq 0$. 

\begin{defn}
We say that $Q_c$ is (nonlinearly) stable in $H^s$ if under \eqref{1p0} one has
\be\label{1p1}
\sup_{t\in\R}\|u(t)-Q_c(\cdot -\rho(t))\|_{H^s} \lesssim \alpha,
\ee
for some $\rho(t)\in\R$.  Otherwise we say that $Q_c$ is unstable.
\end{defn}

\begin{figure}
\begin{center}
\begin{tikzpicture}[
	>=stealth',
	axis/.style={semithick,->},
	coord/.style={dashed, semithick},
	yscale = 1,
	xscale = 1]
	\newcommand{\xmin}{0};
	\newcommand{\xmax}{9};
	\newcommand{\ymin}{0};
	\newcommand{\ymax}{3};
	\newcommand{\ta}{3};
	\newcommand{\fsp}{0.2};
	\draw [axis] (\xmin-\fsp,0) -- (\xmax,0) node [right] {$x$};
	\draw [axis] (0,\ymin-\fsp) -- (0,\ymax) node [below left] {$t$};
	\filldraw[color=light-gray2] (5,1) circle (0.5); 
	\draw [thick,->] (5,1) -- (7,1.1);
	\filldraw[color=light-gray1] (2,1) circle (0.7); 
	\draw [thick,->] (2,1) -- (2.5,1.5);
	\draw [dashed] (5,0) -- (5,1);
	\draw (5,0) node [below] {$x=\rho(t)$};
	\draw (7,1.1) node [above] {$\rho'(t) \sim c$};
	\draw (2,1) node [below] {size $\sim\alpha$};
\end{tikzpicture}
\end{center}
\caption{Faster solitons are darker and more concentrated; speed is commensurate with arrow length. Here a small soliton perturbs the big one, which remains stable up to a modification of the center of mass.}\label{fig:1}
\end{figure}
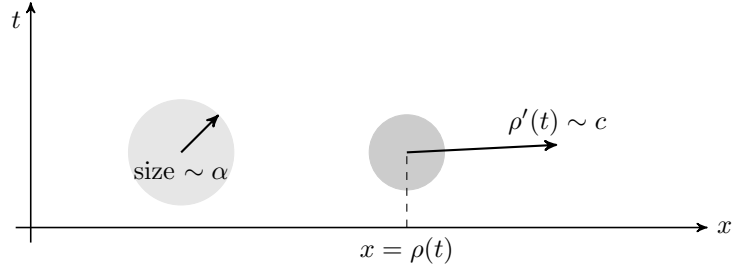

Note that the constant involved in estimate \eqref{1p1} does not depend on $t$ and $\alpha$.  The parameter $\rho(t)$ is absolutely necessary since if $c\sim c'$, with $|c-c'| =\al\ll 1$, one has
\[
\| Q_c-Q_{c'}\|_{H^s} \sim \al \ll1,
\]
but the corresponding solutions satisfy
\[
\| Q_c(\cdot -ct )-Q_{c'}(\cdot -c't)\|_{H^s} \sim 1 ,
\]
as $t\to \infty$.  In that sense, we say that \eqref{1p1} is a sort of \emph{orbital} stability result (see Fig. \ref{fig:1}).

\begin{thm}[$H^1$ stability, \cite{BSS}]\label{Thm1}
Assume $p=2,3,4$, $c>0$, $x_0\in\R$. There are $\alpha_0>0$, $C_0>0$ such that for all $\al\in (0,\al_0)$ and for all $u_0 \in H^1(\R)$ such that 
\[
\|u_0-Q_c(\cdot -x_0)\|_{H^1} <\alpha,
\]
one has 
\[
\sup_{t\in\R}\|u(t)-Q_c(\cdot -\rho(t))\|_{H^s} \leq C_0 \alpha,
\]
for some $\rho(t)$ which satisfies the estimate
\[
\sup_{t\in\R}|\rho'(t)-c| \leq CC_0 \alpha,
\]
for some $C>0$. 
\end{thm}

The previous theorem was proved by Benjamin \cite{Benj} and Bona, Souganidis and Strauss \cite{BSS}. The proof that we present here is based in the approach introduced by M. Weinstein \cite{We2}. 

\subsection{Proof of Theorem \ref{Thm1}} The key idea is the use of conservation laws. Recall that \eqref{gKdV} has at least two conserved quantities: the mass
\be\label{M}
M[u](t) :=\frac12 \int_\R u^2(t,x)dx = M[u_0],
\ee
and the energy
\be\label{E}
E[u](t) :=\frac12 \int_\R u_x^2(t,x)dx -\frac1{p+1}\int_\R u^{p+1}(t,x)dx =E[u_0].
\ee
Both quantities are preserved by the $H^1$ flow. Now we analyze the evolution of a solution $u(t)$ which is a small perturbation of a soliton: if for some $\rho(t)$ fixed we have
\[
u(t,x) = Q_c(x-\rho(t)) +z(t,x),
\]
(note that $z(t,x)$ depends on the definition of $\rho(t)$) then
\bee
E[u](t) & =&  E[Q_c] +\int_\R Q_c' z_x -\int_\R Q_c^p z  \\
& & +\frac12 \int_\R z_x^2 -\frac p2\int_\R Q_c^{p-1}z^2 + O(\|z(t)\|_{H^1(\R)}^3).
\eee
Note that for the first order term on $z$ we have, using \eqref{S},
\[
\int_\R Q_c' z_x -\int_\R Q_c^p z = -c\int_\R Q_c z,
\] 
which can be cancelled if we use the mass
\[
c M[u](t) = cM[Q_c] + c\int_\R Q_c z + \frac c2 \int_\R  z^2. 
\]
Therefore, we have
\bea\label{EM}
E[u_0] + c M[u_0] & =&  E[Q_c] +c M[Q_c]  \nonu\\
& & +\frac12 \int_\R z_x^2 + \frac c2 \int_\R  z^2 -\frac p2\int_\R Q_c^{p-1}z^2 + O(\|z(t)\|_{H^1(\R)}^3).
\eea
Now we consider the quadratic term on $z$  above. It is not difficult to see that 
\[
\frac12 \int_\R z_x^2 + \frac c2 \int_\R  z^2 -\frac p2\int_\R Q_c^{p-1}z^2 = \int_\R z\mathcal Lz,
\] 
where
\be\label{L}
\mathcal Lz := -z_{xx} +cz -pQ_{c}^{p-1}z.
\ee
Now, using classical analysis of operators (see e.g Reed-Simon, vol. 4), we have

\begin{lem}\label{L1p4} The following are satisfied.
\ben
\item $\mathcal L$ is a self-adjoint operator defined on $L^2(\R)$, with domain $D(\mathcal L) =H^2(\R)$. 
\item $\ker \mathcal L $ is spanned by $Q_c'$ (direct from \eqref{S} and ODE analysis). 
\item $\mathcal L$ has a unique negative eigenvalue $-\la_0<0$ of multiplicity one. The associated (explicit) eigenfunction $\chi_c$ satisfies $\chi_c \in S(\R)$ (the Schwartz class), and
\be\label{Lchic}
\mathcal L \chi_c =-\la_0 \chi_c, \quad \|\chi_c\|_{L^2(\R)}=1.
\ee
\item The continuum spectrum of $\mathcal L$ is the closed interval $[c,+\infty)$.
\item \emph{Coercivity}. There is $\ga_0>0$ such that the following holds. Assume that 
\be\label{coer0}
\int_\R \tilde z Q_c' = \int_\R \tilde z \chi_c =0.
\ee
Then 
\be\label{coer1}
\int_\R \tilde z\mathcal L \tilde z \geq \ga_0 \| \tilde z\|_{H^1(\R)}^2.
\ee
\een
\end{lem}

As an immediate corollary of  \eqref{coer1}, one has the following. 

\begin{cor}
There is $\ga_0>0$ such that the following holds. Assume  that
\[
\int_\R \tilde z Q_c'  =0.
\]
Then 
\be\label{coer2}
\int_\R \tilde z\mathcal L \tilde z \geq \ga_0 \| \tilde z\|_{H^1(\R)}^2 -\frac1{\ga_0} \Big| \int_\R \chi_c \tilde z\Big|^2.
\ee
\end{cor}
For the proof of this result, just modify $\tilde z$ with a suitable linear combination of $\chi_c$ such that now \eqref{coer0} is satisfied.

\medskip

However, from the original ideas of Weinstein, we will prove the following very useful variation of \eqref{coer2}.

\begin{lem}\label{L1}
There is $\tilde\ga_0>0$ such that the following holds. Assume  now that
\[
\int_\R z Q_c'  =0.
\]
Then 
\be\label{coer3}
\int_\R z\mathcal Lz \geq \tilde \ga_0 \|z\|_{H^1(\R)}^2 -\frac1{\ga_0} \Big| \int_\R Q_c z\Big|^2.
\ee
\end{lem}

Let us assume for a moment the validity of this lemma. Coming back to the proof of the stability result, we are going to use the Implicit Function Theorem in order to ensure that we can fix $\rho(t)$ such that, for all time $t$,
\be\label{Ortho1}
\int_\R z(t,x) Q_c'(x-\rho(t))dx  = 0,
\ee
so that we fix now $z(t,x)$. Usually we call this method \emph{modulation}. This can be done provided $z(t)$ is small, which is true at least for some time $0<t\leq T_0$. The idea is to prove that we can take $T_0=+\infty$ by a simple bootstrap argument. Note that \eqref{Ortho1} is justified since the smooth functional 
\[
H^1(\R)\times \R \ni (u,\rho) \mapsto  \int_\R(u(x) - Q_c(x-\rho)) Q_c'(x-\rho)dx \in\R
\]
has nondegenerate partial derivative with respect to $\rho$ at the point $(Q_c,\rho)$, because
\[
 \int_\R Q_c'^2  \neq 0.
\]
Therefore, in a small $H^1$-neighborhood of $Q_c$ one has \eqref{Ortho1} well-defined (the uniformity in time for $u(t)\in H^1(\R)$ can be assured at least for a suitable amount of time that will be bootstrapped later on). Using \eqref{coer3} applied to $z(t,x)$ we have
\[
\int_\R z(t) \mathcal Lz(t) \geq \ga_0 \|z(t)\|_{H^1(\R)}^2 -\frac1{\ga_0} \Big| \int_\R Q_c z(t)\Big|^2.
\]
Now we evaluate \eqref{EM} at time $t=0$ and some $t\leq T_0$. We have (we emphasize that the constants involved are independent of time)
\[
 \|z(t)\|_{H^1(\R)}^2 \lesssim  \abs{ \int_\R Q_c z(t) }^2  + \|z(0)\|_{H^1(\R)}^2 +  \|z(t)\|_{H^1(\R)}^3.
\]
Note that if $\al_0$ is chosen smaller, the last term  on the right above can be absorbed by the left hand side. We are left to give a suitable estimate on the first term on the right. But now we use the conservation of mass: we have
\[
\int_\R zQ_c (t) =\int_\R zQ_c (0) + \frac 12 \int_\R z^2(0) -  \frac 12 \int_\R z^2(t), 
\] 
so that 
\[
\abs{\int_\R zQ_c (t) } \lesssim \|z(0)\|_{H^1(\R)} +  \|z(t)\|_{H^1(\R)}^2.
\]
Note that this estimate is better than the usual Cauchy-Schwarz estimate for the term in the left-hand side, since the term of order one in $z(t)$ has disappeared. We finally get
\[
 \|z(t)\|_{H^1(\R)} \lesssim    \|z(0)\|_{H^1(\R)},
\]
so that the estimate on $z(t)$ does not depends on $t$, improving the previous assumption $T_0<+\infty$. This proves the result.

\medskip

Now we prove  Lemma \ref{L1}.

\begin{proof}
We assume that 
\[
\int_\R \tilde z Q_c'  = \int_\R z Q_c = 0.
\]
We want to show that 
\[
\int_\R z \mathcal L z \geq \tilde \ga_0 \|z\|_{H^1(\R)}^2,
\]
for some $\tilde \ga_0>0$ independent of $z$. We decompose $z$ and $\Lambda Q_c$ (cf. \eqref{LaQ}) as follows
\[
z =\bt_0 \chi_c + \tilde z, \quad \bt_0\in \R, \quad  \int_\R \tilde z \chi_c =0,
\] 
and
\[
\Lambda Q_c  =\bt_1 \chi_c + \widetilde{\Lambda Q_c}, \quad \bt_1 \in \R, \quad  \int_\R \widetilde{\Lambda Q_c} \chi_c =0.
\]
Note that we do not need any component in the $Q_c'$ direction since $z$ and $\Lambda Q_c$ are orthogonal to this function.

\medskip

Now, using \eqref{Lchic} and the previous decomposition, we have
\be\label{star0}
\int_\R z\mathcal Lz = -\bt_0^2 \la_0 +\int_\R \tilde z\mathcal L \tilde z.
\ee
Since $p<5$,
\be\label{critC}
\int_\R Q_c \Lambda Q_c>0,
\ee
and
\be\label{star1}
0> -\int_\R Q_c \Lambda Q_c  = \int_\R \Lambda Q_c  \mathcal L \Lambda Q_c  = -\bt_1^2 \la_0 + \int_\R \widetilde{\Lambda Q_c} \mathcal L \widetilde{\Lambda Q_c}.
\ee
Finally, 
\[
0 = \int_\R z Q_c = - \int_\R z \mathcal L \Lambda Q_c = \bt_0\bt_1 \la_0 - \int_\R \widetilde{\Lambda Q_c}  \mathcal L \tilde z ,
\]
so that 
\be\label{star2}
 \bt_0^2\bt_1^2 \la_0^2 =\Big(  \int_\R \widetilde{\Lambda Q_c}  \mathcal L \tilde z \Big)^2.
\ee
From \eqref{star0}, \eqref{star1} and \eqref{star2}, we have
\be\label{star3}
\int_\R z\mathcal Lz  = \int_\R \tilde z\mathcal L \tilde z - \frac{\Big( \displaystyle{ \int_\R \widetilde{\Lambda Q_c}  \mathcal L \tilde z} \Big)^2}{\displaystyle{ \int_\R \widetilde{\Lambda Q_c} \mathcal L \widetilde{\Lambda Q_c} + \int_\R Q_c \Lambda Q_c}}.  
\ee
We get the estimate
\[
\int_\R \tilde z\mathcal L \tilde z \geq \ga_0 \|\tilde z\|_{H^1(\R)}^2 \geq 0.
\]
Moreover, the bilinear operator 
\[
\int_\R v \mathcal L w, \quad v,w \in \{\chi_c, Q_c' \}^\perp
\]
defines an inner product in $L^2(\R)^2$. Therefore, we have a Cauchy-Schwarz inequality:
\[
\abs{\int_\R v \mathcal L w}^2 \leq \Big( \int_\R v \mathcal L v\Big)\Big( w \mathcal L w \Big),
\]
with equality if and only if $v$ is parallel to $w$. Using this information we have
\[
\Big(  \int_\R \widetilde{\Lambda Q_c}  \mathcal L \tilde z \Big)^2 \leq \Big(  \int_\R \widetilde{\Lambda Q_c}  \mathcal L \widetilde{\Lambda Q_c} \Big)\Big(  \int_\R \tilde z  \mathcal L \tilde z \Big),
\]
and now \eqref{star3} becomes
\[
\int_\R z\mathcal Lz \geq  \Big( \int_\R \tilde z\mathcal L \tilde z \Big) \Bigg[ 1 -\frac{ \displaystyle{  \int_\R \widetilde{\Lambda Q_c}  \mathcal L \widetilde{\Lambda Q_c}} }{ \displaystyle{\int_\R \widetilde{\Lambda Q_c} \mathcal L \widetilde{\Lambda Q_c} + \int_\R Q_c \Lambda Q_c}} \Bigg].
\]
Since $ \int_\R Q_c \Lambda Q_c>0$, we have that 
\[
1 -\frac{ \displaystyle{  \int_\R \widetilde{\Lambda Q_c}  \mathcal L \widetilde{\Lambda Q_c}} }{ \displaystyle{\int_\R \widetilde{\Lambda Q_c} \mathcal L \widetilde{\Lambda Q_c} + \int_\R Q_c \Lambda Q_c}} > \eta_0>0,
\]
independent of $z$. We have then
\[
\int_\R z\mathcal Lz \geq  \eta_0 \int_\R \tilde z\mathcal L \tilde z  \geq \eta_0\ga_0\|\tilde z\|_{H^1(\R)}^2 \geq 0.
\]
Finally, from \eqref{star0} and the previous inequality we have 
\bee
 \int_\R \tilde z\mathcal L \tilde z   &=&   \frac12\int_\R \tilde z\mathcal L \tilde z    +\frac12 \int_\R \tilde z\mathcal L \tilde z  \\
 &\geq  &  \frac12\eta_0\ga_0\|\tilde z\|_{H^1(\R)}^2 + \frac12 \bt_0^2 \la_0 \\
 & \geq &  \frac12\eta_0\ga_0\|\tilde z\|_{H^1(\R)}^2 +\frac{\bt_0^2}{C_0} \|\chi_c\|_{H^1(\R)}^2 \\
 & \geq & \tilde\ga_0 \| z\|_{H^1(\R)}^2,
\eee
as desired.

\end{proof}

\medskip

Before finishing this section, some remarks are in order:

\medskip

\begin{rem}
If $p\geq 5$ it is well-known that solitons are \emph{unstable} objects, see Bona-Souganidis-Strauss \cite{BSS} for the case $p>5$, and Martel-Merle \cite{MM} for the more difficult case $p=5$. Compare these results with the critical-supercritical nature of the equation for $p\geq 5$ \eqref{crit}.
\end{rem}

\begin{rem}
One can prove stability using the variational characterization of the soliton as the ground state solution of the elliptic ODE
\[
u'' -c u + u^p =0, \quad u>0, \quad u\in H^1(\R).
\]
See e.g. Cazenave-Lions \cite{CL}, Grillakis-Statah-Strauss \cite{GSS}, etc.
\end{rem}

\begin{rem}
Proving stability in $H^s$, $s\neq 1$ is a very difficult problem. If $s\in (0,1)$, there are polynomial bounds (in time) for the growth of the Sobolev norms (see e.g. \cite{RS}), however, note that the fact that $\|Q_c\|_{H^s} \sim 1$ destroys the utility of such results. If $s=0$ and $p=2$ (namely, we consider the KdV equation with $L^2$ initial data), Merle and Vega \cite{MV} proved $L^2$ stability using the Miura transformation (so the proof is deeply non variational). Their proof has been adapted to several other integrable problems (see Mizumachi-Tzvetkov \cite{MT}, Mu\~noz \cite{Mu2}, Alejo-Mu\~noz-Vega \cite{AMV}, Buckmaster-Koch \cite{BK}, among others).  
\end{rem}

\begin{rem}
One can prove convergence to a soliton at infinity in time, a property called \emph{asymptotic stability}. See the works by Pego and Weinstein \cite{PW}, Martel-Merle \cite{MMnon}, and those related to the nonlinear Schr\"odinger equation, not considered in these notes.
\end{rem}

\begin{rem}
One can also prove stability (and instability) of periodic structures, solutions of \eqref{gKdV}, so \eqref{S1} does not hold in general. See e.g the works by Jaime Angulo \cite{JA} and collaborators.   
\end{rem}

\section{Stability of the sum of $N$-solitons. Multi-solitons}

\medskip

In the previous section we studied the stability problem for a simple soliton $Q_c$. The purpose of this lecture is to consider the problem when several solitons interact themselves, in a weak form. Indeed, fix $p=2,3$ or 4, consider the scaling parameters
\be\label{2p0}
0<c_1<c_2<\cdots<c_N,
\ee
and fixed shift parameters
\be\label{2p1}
x_1^0, x_2^0,\ldots, x_N^0.
\ee
We would like to show that any small perturbation of the sum of $N$ solitons at time $t=0$,
\be\label{2p1a}
\sum_{j=1}^N Q_{c_j}(x-x_j^0),
\ee
where $Q_{c_j}$ is solution of \eqref{S1}, leads to a solution $u(t)$ of \eqref{gKdV} which stays close to the sum of $N$-solitons, of the form
\be\label{2p2}
\sum_{j=1}^N Q_{c_j}(x-\rho_j(t)), \quad \rho_j'(t)\sim c_j
\ee
for all time $t\geq 0.$ before continuing, some remarks are in order.

\medskip

\begin{rem}
Contrary to the previous lecture, the energy space $H^1(\R)$ is not well situed for the stability of more than one soliton, as we will see later. In that sense, any stability result in the energy space for the sum of several solitons is a difficult problem to solve. 
\end{rem}

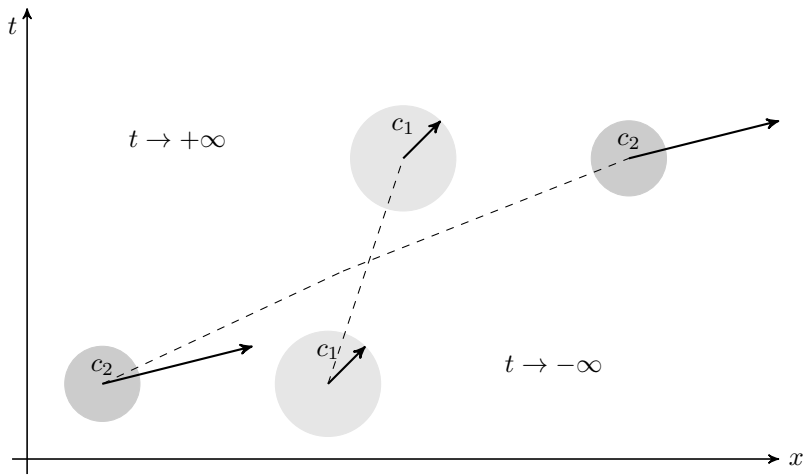
\begin{figure}
\begin{center}
\begin{tikzpicture}[
	>=stealth',
	axis/.style={semithick,->},
	coord/.style={dashed, semithick},
	yscale = 1,
	xscale = 1]
	\newcommand{\xmin}{0};
	\newcommand{\xmax}{10};
	\newcommand{\ymin}{0};
	\newcommand{\ymax}{6};
	\newcommand{\ta}{3};
	\newcommand{\fsp}{0.2};
	\draw [axis] (\xmin-\fsp,0) -- (\xmax,0) node [right] {$x$};
	\draw [axis] (0,\ymin-\fsp) -- (0,\ymax) node [below left] {$t$};
	\filldraw[color=light-gray1] (4,1) circle (0.7); 
	\draw [thick,->] (4,1) -- (4.5,1.5);
	\filldraw[color=light-gray2] (1,1) circle (0.5); 
	\draw [thick,->] (1,1) -- (3,1.5);
	\filldraw[color=light-gray2] (8,4) circle (0.5);
	\draw [thick,->] (8,4) -- (10, 4.5);
	\filldraw[color=light-gray1] (5,4)  circle (0.7); 
	\draw [thick,->] (5,4) -- (5.5,4.5);
	\draw (7,1) node [above] {$t \to -\infty$};
	\draw (2,4) node [above] {$t \to +\infty$};
	\draw (8,4) node [above] {$c_2$};
	\draw (5,4.2) node [above] {$c_1$};
	\draw (4,1.2) node [above] {$c_1$};
	\draw (1,1) node [above] {$c_2$};
	\draw [dashed] (1,1) -- (4.2,2.5);
	\draw [dashed] (4.2,2.5) -- (8,4);
	\draw [dashed] (4,1) -- (5,4);
\end{tikzpicture}
\end{center}
\caption{A 2-soliton for KdV, with $c_1<c_2$. Note that the solution has the same scalings as $t \to \pm \infty.$}\label{fig:2}
\end{figure}

\begin{rem}
If $p=2,3$, it is possible to find \emph{explicit} solutions $U_N(t)$ of \eqref{gKdV} such that 
\be\label{2p3}
\lim_{t\to \pm \infty} \|U_N(t) -\sum_{j=1}^N Q_{c_j}(\cdot -c_j t - x_j^\pm ) \|_{H^1(\R)} =0,
\ee
for $c_j>0$ given and $x_j^\pm \in \R$ explicitly defined (actually, for $p=2,3$ the equation is \emph{completely integrable}, see e.g. \cite{La}). Such solutions are often referred as \emph{multi-solitons}, or $N$-soliton solutions, and describe the interaction (collision) of $N$ different solitons through the evolution in time. It is important to emphasize that, since \eqref{gKdV} is a nonlinear equation, $U(t)$ is always different to the \emph{sum} of $N$ solitons, and \eqref{2p3} holds only at infinity in time. In particular, a good understanding of a property as \eqref{2p2}  is a key step to understand the stability of $U_N(t)$ \emph{globally} in time (see Fig. \ref{fig:2}).
\end{rem}

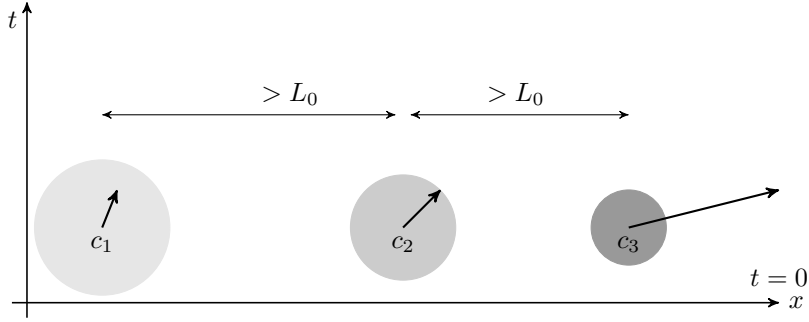
\begin{figure}
\begin{center}
\begin{tikzpicture}[
	>=stealth',
	axis/.style={semithick,->},
	coord/.style={dashed, semithick},
	yscale = 1,
	xscale = 1]
	\newcommand{\xmin}{0};
	\newcommand{\xmax}{10};
	\newcommand{\ymin}{0};
	\newcommand{\ymax}{4};
	\newcommand{\ta}{3};
	\newcommand{\fsp}{0.2};
	\draw [axis] (\xmin-\fsp,0) -- (\xmax,0) node [right] {$x$};
	\draw [axis] (0,\ymin-\fsp) -- (0,\ymax) node [below left] {$t$};
	\filldraw[color=light-gray1] (1,1) circle (0.9); 
	\draw [thick,->] (1,1) -- (1.2,1.5);
	\filldraw[color=light-gray2] (5,1) circle (0.7); 
	\draw [thick,->] (5,1) -- (5.5,1.5);
	\filldraw[color=light-gray3] (8,1) circle (0.5); 
	\draw [thick,->] (8,1) -- (10,1.5);
	\draw (10,0.1) node [above] {$t =0$};
	\draw (5,1) node [below] {$c_2$};
	\draw (1,1) node [below] {$c_1$};
	\draw (8,1) node [below] {$c_3$};
	\draw [<->] (5.1,2.5) -- (8,2.5);
	\draw (6.5,2.5) node [above] {$>L_0$};
	\draw [<->] (1,2.5) -- (4.9,2.5);
	\draw (3.5,2.5) node [above] {$>L_0$};
\end{tikzpicture}
\end{center}
\caption{Well-prepared initial data, where $c_1<c_2<c_3$.}\label{fig:3}
\end{figure}

\begin{rem}
From numerical observations it is clear that a result like \eqref{2p2} cannot hold for the case $p=4$ (the equation is in fact \emph{nonintegrable}), unless the $x_j^0$ in \eqref{2p1} are well-ordered, in the sense that  
\be\label{2p4}
x_1^0< x_2^0<\cdots< x_N^0.
\ee
In fact, under this condition, \emph{multi-collisions} are avoided (see Fig. \ref{fig:3}). 
\end{rem}

The first stability result in the case of multi-solitons was proved by Maddocks and Sachs \cite{MS}, based on an early work by Lax \cite{LAX1}. More precisely, they showed that for $p=2$ (KdV), the multi-soliton $U_N(t)$ in \eqref{2p3} is stable under small perturbations in $H^N(\R)$ (see \cite{MS} for a precise description of the result).  Next, Martel, Merle and Tsai showed that the sum of $N$-solitons as in \eqref{2p1a} is $H^1$ stable. More precisely,

\begin{thm}[\cite{MMT}]\label{Thm2}

Assume that $p=2,3$ or $4$, and that \eqref{2p0} and \eqref{2p4} holds, such that additionally
\be\label{star}
\inf_{j} (x_{j+1}^0 -x_j^0) \geq L_0>0, 
\ee
for some $L_0>0$ large. Then there are $C_0>0$ and $\al_0,\ga_0>0$ such that for all $\al\in (0,\al_0)$, the following holds. Assume that $u_0\in H^1(\R)$ obeys the estimate
\be\label{2p5}
\| u_0 - \sum_{j=1}^N Q_{c_j}(\cdot-x_j^0)\|_{H^1(\R)} < \al,
\ee
then the solution $u(t)$ of \eqref{gKdVCP} with initial data $u_0$ satisfies
\be\label{2p6}
\sup_{t>0} \| u(t) - \sum_{j=1}^N Q_{c_j}(\cdot -\rho_j(t)) \|_{H^1(\R)}  < C_0( \al + e^{-\ga_0 L_0}).
\ee
\end{thm}

Some remarks are in order.

\medskip

\begin{rem}
In the case $p=4$ (quartic gKdV), estimate \eqref{2p6} holds only for $t>0$.  Moreover, condition \eqref{star} is essential. 
\end{rem}

\begin{rem}
Using \eqref{2p3}, estimate \eqref{2p6} improves the Maddocks-Sachs result. Additionally, Theorem  \ref{Thm2} does not use the \emph{integrability} of the equation.  
\end{rem}

\begin{rem}
The previous theorem has been adapted to several dispersive models with soliton solutions, see e.g. Martel-Merle-Tsai \cite{MMT} for the case of the NLS equation, Mu\~noz \cite{Mu2} for the case of the defocusing mKdV equation, among others.
\end{rem}

\begin{rem}
Recently,  Alejo, Mu\~noz and Vega  \cite{AMV} showed that \eqref{2p6} holds even for $L^2$ perturbations, in the case $p=2$ (KdV). In fact, the proof generalizes the Martel-Vega's original result. 
\end{rem}

\subsection{Proof of Theorem \ref{Thm2}} For simplicity, we assume $N=2$. We define
\[
R(t,x) := Q_{c_1}(x-\rho_1(t)) + Q_{c_2}(x-\rho_2(t)), \quad 0<c_1<c_2,
\]
for some $\rho_1(t), \rho_2(t)$ to be chosen later. As in the previous section, the idea is to understand the behavior  of the energy of $u(t)$ in the vicinity of $R(t)$, since $u(0) \sim Q_{c_1}(x-x_1^0) + Q_{c_2}(x-x_2^0)$. 

\medskip

Let us assume that $u(t) =R(t) +z(t)$, for some $z(t)$ small. We have ($z=z(t)$)
\bea\label{2p7}
E[u]& =&  E[R] +\int_\R R_x z_x -\int_\R R^p z \nonu \\
& & +\frac12 \int_\R z_x^2 -\frac p2\int_\R R^{p-1} z^2 +O(\|z\|_{H^1(\R)}^3).
\eea 
We will assume now  that $\rho_1$ and $\rho_2$ are well-separated, at least for a minimum amount of time. Indeed, assume that there is $c_0>0$ such that 
\be\label{rho21}
\rho_2(t) -\rho_1(t) \geq c_0 (t +L_0).
\ee
Using a standard bootstrap argument (in time), this hypothesis will be improved afterwards.  

\medskip

Now we consider the first order term in \eqref{2p7}. We have
\[
\int_\R R_x z_x -\int_\R R^p z =-\int_\R (R_{xx} +R^p)z. 
\] 
Note that from the equation satisfied by each soliton $Q_{c_j}$, we have
\[
R_{xx} +R^p  = -c_1Q_{c_1} -c_2 Q_{c_2}  + [(Q_{c_1} +Q_{c_2})^p -Q_{c_1}^p -Q_{c_2}^p].
\]
Thanks to condition \eqref{rho21} the term between brackets above satisfies the estimate
\[
\|(Q_{c_1} +Q_{c_2})^p -Q_{c_1}^p -Q_{c_2}^p\|_{L^2(\R)} \lesssim e^{-\ga_0 (L_0+t)},
\]
for some $\ga_0>0$. So we have
\[
\int_\R R_x z_x -\int_\R R^p z = -\int_\R (c_1Q_{c_1} + c_2 Q_{c_2})z +O(\|z\|_{L^2(\R)} e^{-\ga_0 (L_0+t)}). 
\]
However, unlike the previous section, the mass cannot help us as usual, since
\bee
M[u](t) & =&  M[Q_{c_1} + Q_{c_2}] + \int_\R (Q_{c_1} + Q_{c_2})z +\frac12 \int_\R z^2\\
& =&  M[Q_{c_1}] + M[Q_{c_2}]+ \int_\R (Q_{c_1} + Q_{c_2})z +\frac12 \int_\R z^2 + O(e^{-\ga_0 (L_0+t)}),
\eee
and $c_1\neq c_2$. The only option  that we have to continue is to impose that the directions associated to each soliton vanish for all time, i.e.
\be\label{Ortho4}
\int_\R Q_{c_1}z = \int_\R Q_{c_2}z =0.
\ee
We would like to use the Implicit Function Theorem to ensure that these orthogonality conditions do hold. A simple computation shows that we cannot modulate in time the variables $\rho_1$ and $\rho_2$ to obtain \eqref{Ortho4} mainly because the direction $Q_{c_j}$ is even and the direction associated to $\rho_j$ is odd $(=Q_{c_j}')$.

\medskip

The idea now is to make $c_1$ and $c_2$ time-dependent functions, i.e. we will modulate the scalings. It is not difficult to check that thanks to \eqref{critC} we can choose $c_1(t)$ and $c_2(t)$ \emph{close to the original scalings} that we denote know $c_1^0$ and $c_2^0$, such that now
\be\label{Ortho4b}
\int_\R Q_{c_1(t)}(x-\rho_1(t)) z (t,x)dx= \int_\R Q_{c_2(t)}( x-\rho_2(t))z(t,x)dx =0,
\ee
and shift parameters $\rho_1(t)$ and $\rho_2(t)$ satisfying
\be\label{Ortho4c}
\int_\R Q_{c_1(t)}'(x-\rho_1(t)) z (t,x)dx= \int_\R Q_{c_2(t)}'( x-\rho_2(t))z(t,x)dx =0,
\ee
where $z(t)$ is given now by
\[
z(t,x) := u(t,x) -\tilde R(t,x), \quad \tilde R(t,x) :=  Q_{c_1(t)}(x-\rho_1(t)) +Q_{c_2(t)}(x-\rho_2(t)).
\]
Note that this choice ensures \eqref{rho21}, at least for certain amount of time.
Coming back to \eqref{2p7}, we have
\bee
E[u_0] & =&  E[\tilde R](t) +\frac12 \int_\R z_x^2(t) -\frac p2 \int_\R \tilde R^{p-1}z^2(t) \\
& & + O(\|z(t)\|_{H^1(\R)}^3 + e^{-\ga_0 (t+L_0)}).
\eee
Note that now $E[\tilde R](t)$ is a function depending on time, so in some sense our problem has became even more difficult to manage, since the O(1) terms in $z$ are now time dependent.

\medskip 

Let us describe more carefully the term $E[\tilde R](t).$ It is not difficult to show that
\[
E[\tilde R](t) = E[Q_{c_1(t)}] + E[Q_{c_2(t)}]  +O(e^{-\ga_0(t+L_0)}),
\]
after using \eqref{rho21} and the fact that $c_1(t)$ and $c_2(t)$ are sufficiently close to the original scalings. Moreover, we have
\be\label{2p11}
E[Q_{c_j(t)}] = c_j^\theta(t) E[Q], \quad \theta := \frac2{p-1}+\frac 12. 
\ee
For the sake of completeness, we also have
\be\label{2p11a}
M[Q_{c_j(t)}] = c_j^{\tilde\theta}(t) E[Q], \quad \tilde\theta := \frac2{p-1}-\frac 12.
\ee
Now the problem will be to control the variation of the quantities $c_j^\theta(t)$ for all time. 

\medskip

Note that 
\bea\label{2p12}
0 = E[u](t) -E[u_0]  & =&  E[\tilde R](t)-E[\tilde R](0)  +\frac12 \int_\R z_x^2(t) -\frac p2 \int_\R \tilde R^{p-1}z^2(t) \nonu \\
& & + O(\|z(0)\|_{H^1(\R)}^2 + \|z(t)\|_{H^1(\R)}^3 + e^{-\ga_0 L_0})
\eea
Now, using \eqref{2p11}
\be\label{2p12a}
E[\tilde R](t)-E[\tilde R](0) = \sum_{j=1,2} (c_j^\theta(t) - c_j^\theta(0)) E[Q] + O(e^{-\ga_0 L_0}).
\ee
From this identity and \eqref{2p12} we have for $j=1,2$,
\be\label{quad0}
|\sum_j \Delta c_j(t)| := | \sum_j c_j(t) -c_j(0)| \lesssim \|z(t)\|_{H^1(\R)}^2 +e^{-\ga_0 L_0},
\ee
so $\sum_j \Delta c_j(t)$ varies \emph{quadratically} on $z(t)$, and in particular it is smaller than $\|z(t)\|_{H^1(\R)}$ (see \cite{MMT} for a proof of the fact that each $\Delta c_j(t)$ is quadratic on $z$).

\medskip

Here comes one of the new ideas in the Martel-Merle-Tsai's paper. Now we estimate \eqref{quad0} using \eqref{2p11a}: We have
\be\label{2p12b}
c_j^\theta(t) - c_j^\theta(0) = \frac{(p+3)}{(5-p)}c_j(0) (c_j^{\tilde\theta}(t) -c_j^{\tilde\theta}(0)) +O(|\Delta c_j(t)|^2),
\ee
(recall that $\theta = \frac2{p-1}+\frac 12$ and $\tilde\theta = \frac2{p-1}-\frac 12$). Now we need an estimate on $c_j^{\tilde\theta}(t) -c_j^{\tilde\theta}(0)$. The second idea in the Martel-Merle-Tsai's paper is to consider a monotonicity property associated to the mass around the sum of $N$ solitons. 

\medskip

\begin{figure}
\begin{center}
\begin{tikzpicture}[
	>=stealth',
	axis/.style={semithick,->},
	coord/.style={dashed, semithick},
	yscale = 1,
	xscale = 1]
	\newcommand{\xmin}{0};
	\newcommand{\xmax}{10};
	\newcommand{\ymin}{0};
	\newcommand{\ymax}{3};
	\newcommand{\ta}{3};
	\newcommand{\fsp}{0.2};
	\draw [axis] (\xmin-\fsp,0) -- (\xmax,0) node [right] {$x$};
	\draw [axis] (0,\ymin-\fsp) -- (0,\ymax) node [below left] {$t$};
	\filldraw[color=light-gray2] (5,1) circle (0.7); 
	\draw [thick,->] (5,1) -- (5.5,1.5);
	\filldraw[color=light-gray3] (8,1) circle (0.5); 
	\draw [thick,->] (8,1) -- (10,1.5);
	\draw [<->] (6.5,-0.5) -- (10,-0.5);
	\draw [<->] (0.2,2.5) -- (10,2.5);
	\draw (8,-0.5) node [below] {$M_2$};
	\draw (4,2.5) node [above] {$M_1=M[u]$};
	\draw (7.3,2.8) node [right] {$x =m_0 t$};
	\draw [dashed] (6.5,-0.6) -- (7.3,3);
	\draw (2,1) node [above] {$t \gg 1$};
	\draw (5,1) node [below] {$c_1$};
	\draw (8,1) node [below] {$c_2$};
\end{tikzpicture}
\end{center}
\caption{The portion of mass $M_2$ for two solitons. Scales are not proportional.}\label{fig:4}
\end{figure}
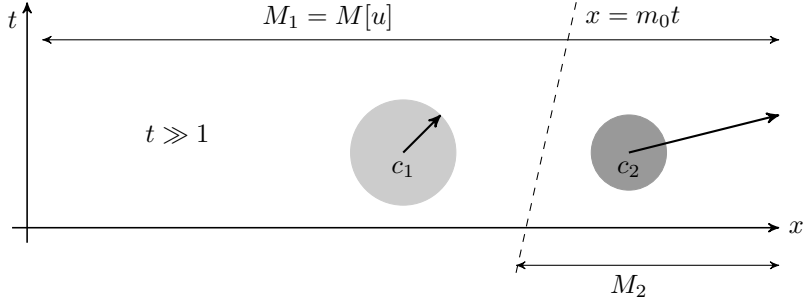

Define
\[
m_0 := \frac12 (c_1^0+ c_2^0),
\] 
and for $A>0$ large but fixed, we consider (see Fig. \ref{fig:4}) 
\[
M_2(t) := \frac 12 \int_\R u^2(t,x) \varphi \Big(\frac{x-m_0 t}{A} \Big),
\]

where $\varphi$ is a smooth kink satisfying
\[
\varphi'>0, \quad  \varphi \in[0,1], \quad  \lim_{-\infty}\varphi =0,\quad \lim_{+\infty}\varphi =1,
\] 
and
\[
|\varphi''' | \leq C\varphi'.
\]
We can take for instance $\varphi(s) = \frac2\pi \arctan e^s$, so that 
\[
\varphi' (s) = \frac{2e^s}{\pi (1+e^{2s})} = \frac{1}{\pi \cosh s}, \quad \varphi''(s) =\frac{- \sinh s}{\pi \cosh^2 s},
\]
and
\[
\varphi'''(s) =\frac{\sinh^2 s-1}{\pi \cosh^3 s} = \frac{1}{\pi \cosh s} . \frac{\sinh^2 s-1}{\cosh^2 s},
\]
so that $|\varphi'''(s)| \leq 2\varphi'(s)$. 

\medskip

Denote 
\[
y:= \frac{x-m_0 t}{A}.
\]
Using the fact that $m_0 t $ represents the \emph{middle point} between both solitons and therefore the \emph{supports} of  $Q_{c_1}(x-\rho_1(t))$ and $\varphi \Big(\frac{x-m_0 t}{A} \Big)$ are largely disjoint, we have
\bea
M_2(t) &  = &  \frac 12 \int_\R (\tilde R +z)^2(t,x) \varphi(y)dx  \nonu\\
& =&  \frac 12 \int_\R (Q_{c_2}^2 +  2Q_{c_2} z + z^2)(t,x) \varphi (y)dx +O(e^{-\ga_0(t+ L_0)})  \nonu\\
& =&  \frac 12 \int_\R (Q_{c_2}^2  + z^2)(t,x) \varphi (y)dx +O(e^{-\ga_0(t+ L_0)}) \nonu \\
& =& c_2^{\tilde\theta}(t) M[Q] +  \frac 12 \int_\R z^2(t,x) \varphi(y) dx +O(e^{-\ga_0(t+ L_0)}), \label{2p13}
\eea
for some $\ga_0>0$ depending on $A$. Note that we have used \eqref{Ortho4} to get rid of the linear term on $z$.

\medskip

On the other hand, following the original Kato smoothing estimate we have
\bea
M_2'(t) &  = &  -\frac{m_0}{2A} \int_\R u^2\varphi' +\int_\R uu_t\varphi \nonu \\
& =& \int_\R (u\varphi)_x (u_{xx} +u^p) - -\frac{m_0}{2A} \int_\R u^2\varphi'   \nonu \\
& =&  -\frac{3}{2A} \int_\R u_x^2 \varphi' +\frac{p}{(p+1)A} \int_\R u^{p+1}\varphi' +\frac1{2A^3}\int_\R u^2\varphi'''  -\frac{m_0}{2A} \int_\R u^2\varphi' .\nonu\\
& &  \label{2p14}
\eea 
Now we estimate \eqref{2p14}. Since $u(t) =\tilde R(t) +z(t)$ we have
\bee
\abs{\int_\R u^{p+1}\varphi'} & \lesssim & \int_\R |\tilde R|^{p+1}\varphi' + \int_\R |z|^{p+1}\varphi'  +e^{-\ga_0(L_0+t)} \\
& \lesssim &  \|z(t)\|_{H^1(\R)}^2 \int_\R z^{2}\varphi'   +  e^{-\ga_0(L_0+t)} .
\eee
Recall that $\varphi' (y)$ and $\tilde R$ have almost disjoint supports and the $L^\infty$ norm is controlled by the $H^1$ norm in one dimension. We also have
\[
\frac1{A^3}\abs{\int_\R u^{2}\varphi''' }\leq \frac2{A^3}  \int_\R u^{2}\varphi' \leq  \frac{m_0}{4A}  \int_\R u^{2}\varphi',
\] 
provided $A$ is large enough. Now \eqref{2p14} can be bounded as follows:
\[
M_2'(t) \leq  -\frac{3}{2A} \int_\R u_x^2 \varphi'   -\frac{m_0}{4A} \int_\R u^2\varphi'   +C \|z(t)\|_{H^1(\R)}^2 \int_\R z^{2}\varphi'  + C e^{-\ga_0(L_0+t)},
\]
so we get
\[
M_2'(t) \leq C  e^{-\ga_0(L_0+t)}.
\]
Integrating in time from $t=0$ and any other time $t>0$, we obtain the monotonicity estimate
\[
M_2(t) \leq M_2(0) +C e^{-\ga_0L_0}.
\]
We use now \eqref{2p13} to get
\[
(c_2^{\tilde\theta}(t) -c_2^{\tilde\theta}(0)) M[Q] +  \frac 12 \int_\R z^2 \varphi (t) \leq C \|z(0)\|_{H^1(\R)}^2 +  C \ e^{-\ga_0L_0}.
\]
Consequently,
\bea\label{C2}
- (c_2(0) -c_1(0))(c_2^{\tilde\theta}(t) -c_2^{\tilde\theta}(0)) M[Q]  &\geq &   \frac 12(c_2(0) -c_1(0)) \int_\R z^2 \varphi (t) \nonu \\
& & - C \|z(0)\|_{H^1(\R)}^2 -  C \ e^{-\ga_0L_0}.
\eea

\medskip

On the other hand, using the conservation of mass we obtain
\bee
& & (c_1^{\tilde\theta}(t) -c_1^{\tilde\theta}(0)) M[Q]  + (c_2^{\tilde\theta}(t) -c_2^{\tilde\theta}(0)) M[Q]  \nonu\\
& & \qquad +  \frac 12 \int_\R z^2 (t) - \frac 12 \int_\R z^2 (0)  +O(e^{-\ga_0 L_0})=0,
\eee
or
\bea\label{2p16a}
& &- c_1(0) (c_1^{\tilde\theta}(t) -c_1^{\tilde\theta}(0)  + c_2^{\tilde\theta}(t) -c_2^{\tilde\theta}(0))M[Q]  \geq  \nonu \\
& &   \qquad \quad  \geq   \frac 12 c_1(0) \int_\R z^2 (t)  - C \|z(0)\|_{H^1(\R)}^2  - C  e^{-\ga_0L_0}.
\eea
Note that 
\bee
 c_1(0) (c_1^{\tilde\theta}(t) -c_1^{\tilde\theta}(0))  +c_2(0) (c_2^{\tilde\theta}(t) -c_2^{\tilde\theta}(0))  &  = & (c_2(0) -c_1(0))(c_2^{\tilde\theta}(t) -c_2^{\tilde\theta}(0)) +  \\
& &+ c_1(0) (c_1^{\tilde\theta}(t) -c_1^{\tilde\theta}(0)  + c_2^{\tilde\theta}(t) -c_2^{\tilde\theta}(0)).
\eee
Now we compute the energy of $Q$. We have from \eqref{S},
\[
-\int_\R Q'^2 =\int_\R QQ'' =\int_\R Q^2 -\int_\R Q^{p+1}.
\] 
Multiplying \eqref{S} by Q'  and integrating, we have
\[
\int_\R Q'^2 = \int_\R Q^2 -\frac{2}{p+1} \int_\R Q^{p+1}.
\]
From these two identities,  we have
\[
\int_\R Q'^2 = \frac{2(p-1)}{p+3}M[Q], \quad \int_\R Q^{p+1} = 4\frac{(p+1)}{p+3} M[Q].
\]
Replacing in $E[Q]$, we obtain
\[
E[Q] = \frac12 \int_\R Q'^2 -\frac1{p+1}\int_\R Q^{p+1} = \frac{p-5}{p+3}M[Q].
\]
Note that the energy of a soliton is negative.  We get the identity 
\[
\frac{(p+3)E[Q]}{(5-p)M[Q]} =-1.
\]
Now \eqref{2p12b} becomes
\[
c_j^\theta(t) - c_j^\theta(0) = \frac{M[Q]}{|E[Q]|}c_j(0) (c_j^{\tilde\theta}(t) -c_j^{\tilde\theta}(0)) +O(\|z(t)\|_{H^1(\R)}^4 + e^{-\ga_0L_0}).
\]
Therefore, using \eqref{2p12a}, \eqref{2p12b}, \eqref{C2} and \eqref{2p16a}
\bee
E[\tilde R](t)-E[\tilde R](0) & =&  -\sum_{j=1,2}c_j(0) (c_j^{\tilde\theta}(t) -c_j^{\tilde\theta}(0)) M[Q]  \\
& & + O(e^{-\ga_0 L_0} + \|z(t)\|_{H^1(\R)}^4) \\
& \geq &  -(c_2(0) -c_1(0))(c_2^{\tilde\theta}(t) -c_2^{\tilde\theta}(0))M[Q] \\
& & - c_1(0) (c_1^{\tilde\theta}(t) -c_1^{\tilde\theta}(0)  + c_2^{\tilde\theta}(t) -c_2^{\tilde\theta}(0))M[Q] \\
& & - Ce^{-\ga_0 L_0} -C \|z(t)\|_{H^1(\R)}^4 \\
& \geq &  \frac 12 c_1(0) \int_\R z^2 (t)(1- \varphi(t)) +\frac 12 c_2(0) \int_\R z^2 (t)\varphi(t) \\
& &  -  C \|z(0)\|_{H^1(\R)}^2 - C \ e^{-\ga_0L_0}.
\eee

Finally, from \eqref{2p12} and the previous estimate we find that
\bea\label{3star}
& & \frac12\int_\R z_x^2(t) +\frac12 c_1(0)\int_\R z^2(1-\varphi)(t) +\frac12 c_2(0)\int_\R z^2 \varphi (t) -\frac p2\int_\R \tilde R^{p-1}z^2 (t)  \nonu \\
& & \qquad \qquad \lesssim \|z(0)\|_{H^1(\R)}^2 +\|z(t)\|_{H^1(\R)}^3   +e^{-\ga_0 L_0}.
\eea
The left hand side of the last inequality can be decomposed in two parts corresponding to each soliton. Indeed, one has
\bee
\tilde R^{p-1}(t) & =&  \tilde R^{p-1}(1-\varphi)(t) + \tilde R^{p-1}\varphi(t)\\
& =& Q_{c_1}^{p-1}(1-\varphi)(t) + Q_{c_2}^{p-1}\varphi(t) + O(e^{-\ga_0 L_0}),
\eee
and
\[
 z_x^2(t) = z_x^2 (1-\varphi)(t) + z_x^2\varphi(t),
\]
so that if $w_1(t) := z(t) \sqrt{1-\varphi (t)}$ and $w_2(t) := z(t) \sqrt{\varphi (t)}$, we get
\bee
z_x^2 \varphi(t) & =&  \Big( -\frac{w_2 \varphi_x }{2\varphi^{3/2}} +\frac{(w_2)_x}{\sqrt{\varphi}} \Big)^2\varphi(t) \\
& =& (w_2)_x^2(t) - \frac{w_2(w_2)_x \varphi_x}{\varphi} + \frac{w_2^2\varphi_x^2}{4\varphi^2}.
\eee
Note that $|\varphi_x/\varphi| \leq C/A$, so for $A$ large we have
\[
\int_\R z_x^2 \varphi(t) = \int_\R  (w_2)_x^2(t)  + O(\frac1A \|w_2(t)\|_{H^1(\R)}^2).
\] 
Similarly,
\[
\int_\R z_x^2 (1-\varphi)(t)  =  \int_\R (w_1)_x^2(t) + O(\frac1A \|w_1(t)\|_{H^1(\R)}^2).
\]
Replacing in \eqref{3star}, we obtain
\bee
& &  \frac12\int_\R (w_1)_x^2(t) +\frac12 c_1(0)\int_\R w_1^2(t) -\frac p2\int_\R Q_{c_1}^{p-1}w_1^2 (t) \\
& & +\frac12\int_\R (w_2)_x^2(t) +\frac12 c_2(0)\int_\R w_2^2(t) -\frac p2\int_\R Q_{c_2}^{p-1}w_2^2 (t) \\
& & \qquad \qquad \lesssim \|z(0)\|_{H^1(\R)}^2 +\|z(t)\|_{H^1(\R)}^3   +e^{-\ga_0 L_0}.
\eee
Finally, we use the orthogonality conditions \eqref{Ortho4b}-\eqref{Ortho4c} as follows: note that
\bee
0& =&  \int_\R Q_{c_2}z \\
& =&   \int_\R Q_{c_2} [ z\sqrt{\varphi}  +z (1-\sqrt{\varphi})] \\
& =&   \int_\R Q_{c_2} w_2 + O(\|z(t)\|_{L^2(\R)}e^{-\ga_0 L_0}),
\eee
with a constant depending on $A$ (we can take $L_0$ larger if necessary). So $w_2$ is almost orthogonal to $Q_{c_2}$, and the error is small enough to be controlled using the formulation of coercivity described in \eqref{coer3}. A similar argument allows to deal with every orthogonality condition. We finally get
\[
\|w_1(t)\|_{H^1(\R)}^2 + \|w_2(t)\|_{H^1(\R)}^2 \lesssim \|z(0)\|_{H^1(\R)}^2 +\|z(t)\|_{H^1(\R)}^3   +e^{-\ga_0 L_0}.
\]
Finally, note that for $A$ large,
\[
\|w_1(t)\|_{H^1(\R)}^2 + \|w_2(t)\|_{H^1(\R)}^2 \sim \|z(t)\|_{H^1}^2.
\]
Using this equivalence, we conclude.

\bigskip

\section{The collision problem for nonintegrable gKdV equations}

\medskip

In the previous Section (see \eqref{2p3}) we introduced the notion of $N$-soliton, or multi-soliton. The fact that the solution decomposes at infinity into the \emph{same} original solitons (with different shifts only) is usually denoted as \emph{elasticity}, that is, the collision among $N$ solitons is elastic.

\medskip

Recall that such a phenomenon is valid only for $p=2$ and 3, although the case $p=3$ is a little bit more complicated to describe, as we will see in the next section.  It is believed that no $N$-soliton solution exists for $p=4$, mainly because the quartic gKdV equation is nonintegrable.

\medskip
 
In 2005, Martel \cite{Martel} proved the existence and uniqueness of \emph{N-soliton-like} solutions for gKdV, $p=2,3$ and 4. More precisely, there is a unique solution $U(t)$ of \eqref{gKdV} such that 
\be\label{3p1}
\lim_{t\to - \infty} \|U(t) -\sum_{j=1}^N Q_{c_j}(\cdot -c_j t - x_j^- ) \|_{H^1(\R)} =0,
\ee
for $(c_j)$ and $x_j^-$ given.  Clearly $U(t)$ coincides with the $N$-soliton solution $U_N(t)$ given in \eqref{2p3}, for $p=2$ and 3. However, for $p=4$ the behavior of $U(t)$ as $t\to +\infty$ is unknown, mainly because of the existence of a strong  regime of interaction. 

\medskip

Martel's idea has shown to be a very robust technique to show existence of $N$-soliton-like solutions for a wide spectrum of dispersive models, mainly because the existence of such objects follows from compactness ideas, and it is not related to the stability of several solitons. In fact, it is possible to construct such solutions even if the corresponding single solitons are unstable! 
See e.g. \cite{CMM} and \cite{CM}.

\medskip

Let us explain the ideas behind \eqref{3p1}.

\medskip

The proof is very similar to the stability proof described in the previous section. However, this time we do not need a monotonicity property since we are just constructing a particular family of solutions, and not describing the behavior of an open set of initial data. 

\medskip

From the fact that the interactions between different speed solitons decrease exponentially in time (due to the exponential decay of the soliton solution), Martel proved the uniform bound in time
\be\label{3p3}
\|U(t) -\sum_{j=1}^N Q_{c_j}(\cdot -c_j t - x_j^- ) \|_{H^1(\R)} \leq C e^{\ga_0 t}, \quad t< 0,
\ee
for some $\ga_0>0$ depending on the \emph{minimal scaling}. Of course the bound above is interesting when $t$ is largely negative, otherwise it loses its effectiveness.

\medskip

In order to prove \eqref{3p3}, the idea is as follows: we consider a decreasing  sequence of times $T_n$ approaching $-\infty$, and we solve the Cauchy problem for \eqref{gKdVCP} with a particular set of initial data:
\[
u_n(T_n) = R(T_n), 
\]
where
\[
R(t,x) := \sum_{j=1}^N Q_{c_j}(x-c_j t - x_j^- ). 
\]
Now the idea is to prove uniform estimates on $n$, for all time relatively large. More precisely, one has
\be\label{3p3a}
\|u_n(t) - R(t)\|_{H^1(\R)} \leq Ce^{\ga_0 t},
\ee
with constants independent of $n$, and for all $t\leq -T_0<0$ fixed. Note that this estimate is indeed  valid for a certain amount of time near each $T_n$. Using the same idea as in the stability proof, one can bootstrap this estimate (without needing the monotonicity property at full) to conclude that  \eqref{3p3a} holds uniformly in time and the bounds are uniform on $n$. Finally, passing to the limit $n\to +\infty$ and using the continuity of the gKdV flow we get the desired estimate.

\medskip

Once \eqref{3p1} is proved, the idea is to study the \emph{collision} problem.  Previous results in that direction can be found in Mizumachi \cite{Mo}, and some numerical simulations suggest that for $p=4$, no elastic collision happens. In 2007, Martel and Merle considered a particular (simpler, but by no means less difficult) problem. They assumed $p=4$, $N=2$ (i.e. only two solitons), and additionally one soliton has to be smaller compared with the other one. Under this situation, they showed for the first time that the collision is indeed \emph{inelastic}. They also consider \cite{MMcol2} the case of general nonlinearities $f(u)$, proving stability of the two-soliton structure, but leaving open the question of inelasticity. Finally, we showed in \cite{Mu1} that no matter what the nonlinearity is (except for the integrable cases), the collision between two small solitons, one being even smaller compared with the other, is always inelastic.

\medskip

\begin{thm}[\cite{MMcol1,MMcol2,Mu1}]\label{Thm3}
Assume $p=4$ (quartic gKdV), $c_2=1$ and $c_1=c\ll1$. Consider the unique solution $U(t)$ satisfying
\[
\lim_{t\to -\infty}\|U(t) - Q(\cdot -t)-Q_c(\cdot -ct)\|_{H^1(\R)} =0.
\] 
Then there are $C_0>0,$ $T_c \gg 1$, $c_1^+>1$ and $c_2^+ \in (0,c)$ (with $c_1^+ \sim 1$ and $c_2^+ \sim c$ in terms of powers of $c$ smaller than $c^{11/12}$), and such that 
\ben
\item Stability. 
\be\label{3p4}
\sup_{t\geq T_c} \|U(t) - Q_{c_1^+}(\cdot - \rho_1(t))-Q_{c_2^+}(\cdot -\rho_2(t))\|_{H^1(\R)} \leq C_0 c^{11/12},
\ee
for some $\rho_1(t)$ and $\rho_2(t)$ satisfying the standard estimates.
\item Nonexistence.
\be\label{3p5}
\liminf_{t\to +\infty} \|U(t) - Q_{c_1^+}(\cdot - \rho_1(t))-Q_{c_2^+}(\cdot -\rho_2(t))\|_{H^1(\R)} \geq \frac1{C_0} c^{17/12}.
\ee
\een
Moreover, a similar result holds for any nonlinearity $f(u)$ in \eqref{gKdV}, provided 
\ben
\item $f$ is subcritical around $u=0$,
\item $f(u) \neq u^2$, $u^3$ and $u^2 +m u^3$, $m\in \R$, and
\item $f$ has solitons $Q_{c_1}$, $Q_{c_2}$ with $0<c_2\ll c_1\ll1$.
\een
\end{thm}

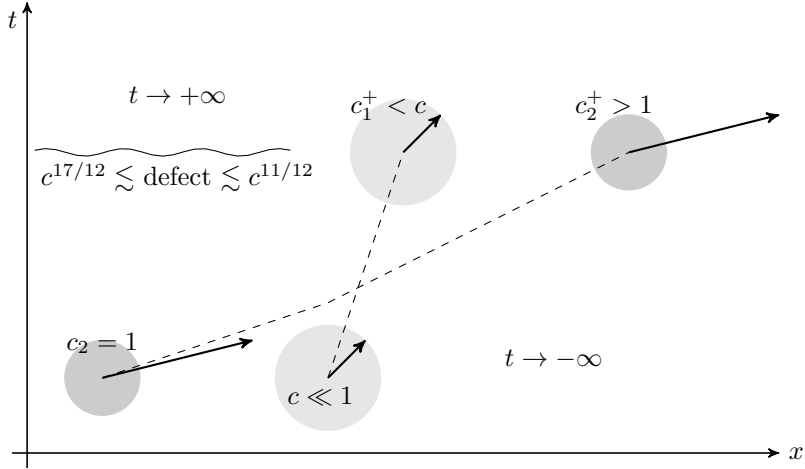
\begin{figure}
\begin{center}
\begin{tikzpicture}[
	>=stealth',
	axis/.style={semithick,->},
	coord/.style={dashed, semithick},
	yscale = 1,
	xscale = 1]
	\newcommand{\xmin}{0};
	\newcommand{\xmax}{10};
	\newcommand{\ymin}{0};
	\newcommand{\ymax}{6};
	\newcommand{\ta}{3};
	\newcommand{\fsp}{0.2};
	\draw [axis] (\xmin-\fsp,0) -- (\xmax,0) node [right] {$x$};
	\draw [axis] (0,\ymin-\fsp) -- (0,\ymax) node [below left] {$t$};
	\filldraw[color=light-gray1] (4,1) circle (0.7); 
	\draw [thick,->] (4,1) -- (4.5,1.5);
	\filldraw[color=light-gray2] (1,1) circle (0.5); 
	\draw [thick,->] (1,1) -- (3,1.5);
	\filldraw[color=light-gray2] (8,4) circle (0.5);
	\draw [thick,->] (8,4) -- (10, 4.5);
	\filldraw[color=light-gray1] (5,4)  circle (0.7); 
	\draw [thick,->] (5,4) -- (5.5,4.5);
	\draw (7,1) node [above] {$t \to -\infty$};
	\draw (2,4.5) node [above] {$t \to +\infty$};
	\draw [domain=0.1:3.5] plot (\x, {4+ sin(350*\x )/20});
	\draw (2,4) node [below] {$c^{17/12} \lesssim$ defect $\lesssim c^{11/12}$};
	\draw (7.8,4.3) node [above] {$c_2^+>1$};
	\draw (4.8,4.3) node [above] {$c_1^+<c$};
	\draw (3.9,1) node [below] {$c\ll 1$};
	\draw (1,1.2) node [above] {$c_2=1$};
	\draw [dashed] (1,1) -- (4,2);
	\draw [dashed] (4,2) -- (8,4);
	\draw [dashed] (4,1) -- (5,4);
\end{tikzpicture}
\end{center}
\caption{Theorem \ref{Thm3} for 2-soliton-like solution of the quartic gKdV equation, with $0<c\ll 1$. Note that the solution has not the same scalings as $t \to \pm \infty.$ The dashed lines are not necessarily straight lines. The defect is deeply related to the differences between final and original scalings.}\label{fig:5}
\end{figure}

Some remarks are in order. 

\begin{rem}
The scaling $c_1^+$ and $c_2^+$ are unique in the sense of the \emph{asymptotic stability} result:
\be\label{TT}
\lim_{t\to +\infty} \|U(t) - Q_{c_1^+}(\cdot - \rho_1(t))-Q_{c_2^+}(\cdot -\rho_2(t))\|_{H^1(x> \frac 12ct)} =0.
\ee
If this condition is lifted, then $c_2^+$ and $c_1^+$ in \eqref{3p4} are not necessarily unique. Note that \eqref{TT} says that in the region $x>\frac12ct$ we have nothing except two solitons, at infinity.
\end{rem}

\begin{rem}
In \cite{MMInv}, Martel and Merle consider the case of collision between nearly equal solitons, proving similar results. In particular, the collision between two solitons is stable, but inelastic. For a similar account on this work, see e.g. \cite{MMEd}. 
\end{rem}

\begin{rem}
The bounds $c^{11/12}$ and $c^{17/12}$ are inherent to the quartic case, in the general case where $u^4$ is replaced by $f(u)$ they may change. However, the difference between them (of the order $c^{1/2}$) seems to be always present. The understanding or improvement of this difference is a nice open problem. 
\end{rem}

The proof of Theorem \ref{Thm3} is involved and very long. The purpose of these notes is to give a suitable account of the main ideas of the proof of \eqref{3p4} and \eqref{3p5}.

\subsection{Sketch of proof of Theorem \ref{Thm3}} Recall that from \eqref{3p3} we have for all $t<0$,
\be\label{3p6}
\|U(t) - Q(\cdot -t) -Q_c(\cdot -ct)\|_{H^1(\R)} \leq C_0 e^{\ga_0 \sqrt{c}t}.
\ee
The constant $\sqrt{c}$ appears in the exponential term after a careful analysis of the interaction between two solitons of size 1 and $c$ in \eqref{3p3}. 
From \eqref{3p6} it is clear that we will have a good control on the solution $U(t)$ provided $t \ll c^{-1/2}$. For this reason we define
\[
T_c := c^{-1/2-\delta_0},
\]
for some $\delta_0>0$ small but fixed. Note that we have
\bea
\|U(-T_c) - Q(\cdot +T_c) -Q_c(\cdot + cT_c)\|_{H^1(\R)} & \leq&  C_0 e^{-\ga_0 \sqrt{c}T_c} \nonu\\
& =& C_0 e^{-\ga_0 c^{-\delta_0}} \ll c^{1000}, \label{mTe}
\eea
for $c\ll1$ small. In other words, the solution is almost the sum of two solitons at time $t=-T_c$. We denote by
\[
[-T_c,T_c]
\]
the interval of interaction in time (note that we do not know if the interaction lasts such an amount of time). We hope to have that at $t=T_c$,
\be\label{3p7}
\|U(T_c) - Q(\cdot -a_+) -Q_c(\cdot -b_+)\|_{H^1(\R)}  \lesssim c^{11/12},
\ee
for some $a_+,b_+ \in \R$ satisfying $a_+ -b_+ \gg \frac 12 T_c$.  Note that estimate \eqref{3p7} makes sense since the soliton satisfies the estimate
\[
\|Q_c\|_{H^1(\R)} \sim c^{1/12},
\]
namely, it is very much larger than the bound $c^{11/12}$.  Assuming \eqref{3p7} and invoking Theorem \ref{Thm2} with $c$ small, we will have \eqref{3p4}.

\medskip

Let us see how the proof of  \eqref{3p7} works.

\medskip

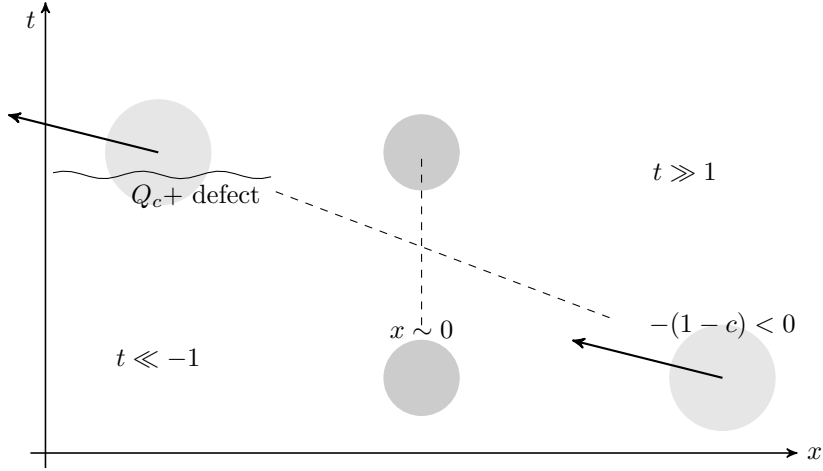
\begin{figure}
\begin{center}
\begin{tikzpicture}[
	>=stealth',
	axis/.style={semithick,->},
	coord/.style={dashed, semithick},
	yscale = 1,
	xscale = 1]
	\newcommand{\xmin}{0};
	\newcommand{\xmax}{10};
	\newcommand{\ymin}{0};
	\newcommand{\ymax}{6};
	\newcommand{\ta}{3};
	\newcommand{\fsp}{0.2};
	\draw [axis] (\xmin-\fsp,0) -- (\xmax,0) node [right] {$x$};
	\draw [axis] (0,\ymin-\fsp) -- (0,\ymax) node [below left] {$t$};
	\filldraw[color=light-gray1] (9,1) circle (0.7); 
	\draw [thick,->] (9,1) -- (7,1.5);
	\filldraw[color=light-gray2] (5,1) circle (0.5); 
	\filldraw[color=light-gray2] (5,4) circle (0.5); 
	\filldraw[color=light-gray1] (1.5,4)  circle (0.7); 
	\draw [thick,->] (1.5,4) -- (-0.5, 4.5);
	\draw (1.5,1) node [above] {$t \ll -1$};
	\draw (8.5,3.5) node [above] {$t \gg 1$};
	\draw [domain=0.1:3] plot (\x, {3.7+ sin(350*\x )/20});
	\draw (2,3.7) node [below] {$Q_c + $ defect};
	\draw (9,1.4) node [above] {$-(1-c)<0$};
	\draw (5,1.4) node [above] {$x\sim 0$};
	\draw [dashed] (5,1.7) -- (5,4);
	\draw [dashed] (7.5,1.8) -- (3,3.5);
\end{tikzpicture}
\end{center}
\caption{Theorem \ref{Thm3} for the function $v(t,x)$}\label{6}
\end{figure}

First of all, we will place ourselves at the origin in space. Let us define 
\[
v(t,x) := u(t,x+t).
\]
It is not difficult to check that $v$ satisfies the equation
\be\label{Sv}
S[v]:= v_t + (v_{xx}-v+v^4)_x =0.
\ee
The collision problem for $v$ is composed by two solitons, $Q(x)$ fixed at the origin and a small soliton with negative velocity $Q_c(x+(1-c)t)$ coming from positive infinite in space (see Fig. \ref{6}). Note that $(1-c)>0$, so the collision is in fact happening \emph{fast} in time, but \emph{slowly} in space, since the essential support of $Q_c$ is  of order $c^{-1/2}$. This last fact also justifies the choice of $T_c$ as the time of interaction.

\medskip

Let us introduce some important notation. We denote
\be\label{y}
y_c := x+ (1-c)t, \quad y:=x-\al(y_c),
\ee
where $\al$ is a function to be chosen later. Recall that 
\[
S[v]:= v_t + (v_{xx}-v+v^4)_x =0.
\]
Now we announce some easy-to-verify but very useful facts. Below $\ga_0$ is a fixed small, positive constant, independent of $c$. We have
\ben
\item For $p=4$,  $Q_c(s) =c^{1/3}Q(\sqrt{c} s)$, and $Q_c''-cQ_c +Q_c^4 =0$ is the equation of the soliton.
\item For any $k\geq 1$ integer, we have
\be\label{solL2}
\|Q_c^k\|_{H^1(\R)} \sim c^{\frac k3 -\frac14}, \quad \|(Q_c^k)'\|_{H^1(\R)} \sim c^{\frac k3 +\frac14}.
\ee
\item Assume now that $f(x) \in \mathcal S(\R)$, then
\be\label{5a}
\|f(x)Q_c^k(y_c)\|_{H^1(\R)} \sim c^{\frac k3}e^{-\ga_0\sqrt{c}|t|},
\ee
and
\be\label{5b}
\|f(x)(Q_c^k)'(y_c)\|_{H^1(\R)} \sim c^{\frac k3 +\frac12}e^{-\ga_0\sqrt{c}|t|}.
\ee
\een

\medskip

Now we propose a first ansatz, say $v_0$, for an approximate solution. Assume 
\[
v_0(t,x) := Q(x) + Q_c(y_c).
\]
Then, using \eqref{S},
\bee
S[v_0] & = & ((Q+Q_c)^4 -Q^4-Q_c^4)_x \\
& =& (4Q^3)'Q_c +4Q^3Q_c' + (6Q^2)' Q_c^2 + 6Q^2(Q_c^2)' +4Q' Q_c^3 +4Q(Q_c^3)'.
\eee
Note that from \eqref{5a} we have that the worst  term is actually the first one, since
\be\label{5c}
\|(4Q^3)'Q_c \|_{H^1(\R)} \lesssim c^{1/3}e^{-\ga_0\sqrt{c}|t|}. 
\ee
In what follows, we will make use of the following \emph{principle of stability}:

\begin{Cl}\label{SP} If $\|S[\tilde v](t)\|_{H^1(\R)}\lesssim c^m$ on $[-T_c,T_c]$, then the actual solution  $v(t)$ should satisfy
\[
\|v(t) -\tilde v(t)\|_{H^1(\R)} \lesssim c^m T_c,
\]
modulo a modulation in time on $v(t)$, and for all $t\in [-T_c, T_c]$. 
\end{Cl}
\medskip

This estimate can be seen as the effect of the propagation of the error term $S[v](t)$ through the interval of time $[-T_c, T_c]$. Note that this principle is in concordance with the case of a standard  soliton $Q(x)$ which is stable for all time, since one has $S[Q]\equiv 0$.

\medskip

For a moment we will assume the validity of this \emph{pseudo theorem}, that we will prove later.  Using this principle and \eqref{5c}, we should have
\[
\|v(t) -\tilde v(t)\|_{H^1(\R)} \lesssim c^{-1/12-\delta_0},
\]
which is extraordinary large. 

\medskip

Now the goal is to get rid of the term $(4Q^3)'Q_c$. Using the algebra related to this term we will assume that 
\[
v_1(t,x) := Q(x) + Q_c(y_c) + A(x) Q_c(y_c),
\]
with $A$ unknown. Replacing in $S[\cdot]$ we obtain now
\bee
S[v_1] & =& (1-c) AQ_c' +[ (AQ_c)_{xx} -AQ_c +4Q^3 AQ_c]_x \\
& & + [ (Q+Q_c+AQ_c)^4 -Q^4-Q_c^4-4Q^3AQ_c]_x .
\eee
Simplifying we obtain
\bea\label{Sv1}
S[v_1] & =&  Q_c(A''-A+4Q^3A +4Q^3)'  \nonu\\
& & + Q_c'(3A''+4Q^3A -cA) + \hbox{smaller terms}.
\eea
Note that $A''-A+4Q^3A +4Q^3 =0$ can be written as (cf. \eqref{L} and Lemma \ref{L1p4})
\[
(\mathcal L A)' = (4Q^3)', \quad \mathcal L :=\mathcal L_{c=1},
\]
so that we have a solution $A \in L^2(\R)$ provided 
\[
\int_\R 4Q^3 Q' =0,
\]
which is indeed the case. Moreover, it is clear that $A$ is in the Schwartz class, and it has exponential decay. More precisely, $A$ is explicitly given by the quantity
\be\label{AA}
A= \frac13 Q'\int_0^x Q^2 -\frac23Q^3.
\ee
In other words (see \eqref{5a}),
\[
\|A(x) Q_c(y_c)\|_{H^1(\R)}  \lesssim  c^{1/3}e^{-\ga_0\sqrt{c}|t|},
\]
which implies that at time $t=T_c$ this term almost disappears, in other words, it cannot represent a defect appeared after the interaction.

\medskip

Now we consider an improvement of $v_1$, such that the second term in \eqref{Sv1} disappears. Following the same idea as before, we consider
\[
v_2(t,x) = Q(x) + Q_c(y_c) + A(x) Q_c(y_c) + B(x) Q_c'(y_c),
\]
with $B$ unknown. Replacing in the equation, we obtain the following equation for $B:$
\[
(\mathcal L B)' = 3A'' +4Q^3 A + 4Q^3 \in \mathcal S(\R).
\]
Note that not every term on the right above is the derivative of a localized function. Even worse, in order to have a solution for the previous linear equation we need
\[
\int_\R (3A'' +4Q^3 A + 4Q^3)Q=0,
\]
since
\[
\int_\R Q(\mathcal L B)' =-\int_\R Q' \mathcal L B = -\int_\R B \mathcal L Q' =0.
\]
However, we have
\[
\int_\R (3A'' +4Q^3 A + 4Q^3)Q \neq 0.
\]
Indeed, note that from \eqref{S} and \eqref{AA},
\bee
\int_\R (3A'' +4Q^3 A + 4Q^3)Q & =&  \int_\R A(3Q'' +4Q^4) +4\int_\R Q^4 \\
& =&\int_\R A(3Q +Q^4) +4\int_\R Q^4 \\
& =&  2\int_\R Q^4 -\frac23 \int_\R Q^7 + \frac13\int_\R (3Q +Q^4) Q' \int_0^x Q^2 \\
& =& 2\int_\R Q^4 -\frac23 \int_\R Q^7 - \frac13\int_\R (\frac 32Q^2 + \frac15Q^5) Q^2 \\
& =& \frac32 \int_\R Q^4 -\frac{11}{15}\int_\R Q^7.
\eee
Finally, we have from \eqref{S},
\[
\int_\R Q^4 = \int_\R Q, \quad \int_\R Q^7 = \frac{20}{11}\int_\R Q,
\]
so that 
\[
\int_\R (3A'' +4Q^3 A + 4Q^3)Q = \frac16\int_\R Q.
\]

\medskip

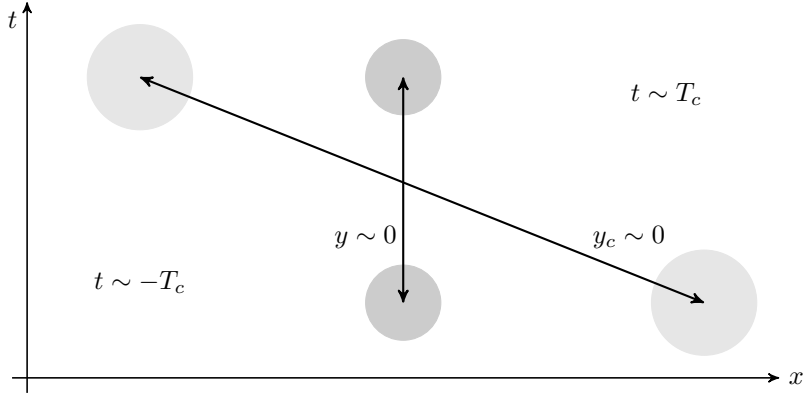
\begin{figure}
\begin{center}
\begin{tikzpicture}[
	>=stealth',
	axis/.style={semithick,->},
	coord/.style={dashed, semithick},
	yscale = 1,
	xscale = 1]
	\newcommand{\xmin}{0};
	\newcommand{\xmax}{10};
	\newcommand{\ymin}{0};
	\newcommand{\ymax}{5};
	\newcommand{\ta}{3};
	\newcommand{\fsp}{0.2};
	\draw [axis] (\xmin-\fsp,0) -- (\xmax,0) node [right] {$x$};
	\draw [axis] (0,\ymin-\fsp) -- (0,\ymax) node [below left] {$t$};
	\filldraw[color=light-gray1] (9,1) circle (0.7); 
	\filldraw[color=light-gray2] (5,1) circle (0.5); 
	\filldraw[color=light-gray2] (5,4) circle (0.5); 
	\filldraw[color=light-gray1] (1.5,4)  circle (0.7); 
	\draw [thick,<->] (9,1) -- (1.5, 4);
	\draw [thick,<->] (5,1) -- (5,4);
	\draw (1.5,1) node [above] {$t\sim -T_c$};
	\draw (8.5,3.5) node [above] {$t \sim T_c$};
	\draw (8,1.6) node [above] {$y_c\sim 0$};
	\draw (4.5,1.6) node [above] {$y\sim 0$};
\end{tikzpicture}
\end{center}
\caption{The two main variables for the collision problem: $y=x-\al(y_c)$ and $y_c= x+(1-c)t$.}\label{9}
\end{figure}

\emph{Key observation:} We forgot the shift on $Q(x)$! Indeed, even in the integrable cases $Q$ has a nontrivial shift.  The idea is to introduce the variable $y$ defined in \eqref{y}, with a function $\alpha$ representing a shift (so it must vary from one nontrivial quantity to another one during a large period of time). In order to preserve the algebra already introduced, the key point is to take $\al(y_c)$ of the form
\[
\al(y_c) = a_1 \int_0^{y_c} Q_c(s)ds. 
\] 
With this choice any derivative of $Q(y)$ will give rise to a new term which is a suitable dilation of $Q_c(y_c)Q'(y)$, and it will preserve the multiplicative algebra. Similarly,
\[
\partial_x Q(y) = Q'(y) -\al'(y_c) Q'(y) = Q'(y) -a_1 Q'(y) Q_c(y_c).
\]
Note that the second term  in the last equality is smaller compared with the first one, so the contribution of the shift will be at the level of the next linear problem. A similar reasoning is valid for $\partial_t Q(y) .$  See Fig. \ref{9} for the meaning of $y$ and $y_c$.

\medskip

We define
\[
v_3(t,x) :=  Q(y) + Q_c(y_c) + A(y) Q_c(y_c) + B(y) Q_c'(y_c).
\]
Now we look for $B$ and $a_1\in \R$. After a standard  computation we obtain the perturbed linear system
\[
(\mathcal L A)' + a_1(3Q-2Q^4)' =(4Q^3)',
\]
and
\[
(\mathcal L B)' + a_1(3Q)'' -3A'' -4Q^3 A =4Q^3.
\]
Recall that making $a_1=0$ we recover the original system, impossible to solve. However, now it is perfectly possible since
\[
\int_\R Q'' Q = -\int_\R Q'^2 \neq 0.
\]
However, the fact that the term $ \mathcal L B$ has no localized  right hand term implies that $B$ must be chosen non localized. In fact, it is possible to solve the previous system as follows:

\begin{lem}
We have
\[
a_1 =-2\frac{\int_\R Q}{\int_\R Q^2}, \quad B=-2\varphi, \quad \varphi := -\frac{Q'}{Q},
\]
and
\[
A =\frac13 Q'\int_0^x Q^2 -\frac23Q^3 -a_1(\frac13Q -\frac32 yQ').
\]
\end{lem}

Note that $B$ is a \emph{kink}, namely an odd function with negative derivative and limits at infinity. The uniqueness of $A$, $B$ and $a_1$ has to be understood in the following sense: $A$ is the unique even, localized solution, and $B$ is odd. 

\medskip

Let us discuss the interpretation of $B$. Note that 
\[
B=-2\varphi(y) \to \mp 2 \quad \hbox{ as } y\to \pm \infty,
\]
so at time $t=-T_c$ we have that if $y_c \sim 0$ then $y \gg 1$, so  
\bee
v_3(-T_c,x) & \sim &  Q(y) + Q_c(y_c) - 2Q_c'(y_c) \\
& \sim & Q(y) + Q_c(y_c -2).
\eee
Similarly, at $t=T_c$ we have $y\ll -1$ if $y_c \sim 0$, and
\bee
v_3(T_c,x) & \sim &  Q(y) + Q_c(y_c) + 2Q_c'(y_c) \\
& \sim & Q(y) + Q_c(y_c + 2).
\eee
In other words, $B$ represents the first order expansion of the shift appearing on the small soliton through the interaction.

\medskip

The problem now is that we have not found an actual defect appearing from the interaction (something not related with shifts, scalings, etc.)

\medskip

We perform a new ansatz. Now we take
\bee
v_4(t,x) &:=&  Q(y) + Q_c(y_c) +A_1(y) Q_{c}(y_c) +B_1(y) Q_{c}'(y_c) \\
& &   +A_2(y) Q_{c}^2(y_c) +B_2(y) (Q_{c}^2)'(y_c),
\eee
where $A_1=A$ and $B_1 =B$ as before, $A_2$ and $B_2$ are unknown functions, and 
\be\label{YAL}
y:= x-\al(y_c), \quad \al(y_c) := a_1\int_0^{y_c} Q_c(s)ds + a_2\int_0^{y_c} Q_c^2(s)ds.
\ee
After replacing in the equation, we will obtain the following linear system for $A_2$, $B_2$ and $a_2$:
\be\label{Omega2}
\begin{cases}
(\mathcal L A_2)' +a_2(3Q-2Q^4)'  =F_2,\\
(\mathcal L B_2)'  +3a_2 Q'' -3A_2'' -4Q^3A_2  =G_2.
\end{cases}
\ee 
Here 
\[
F_2 := (6Q^2 (1+A_1)^2)' -a_1(4Q^3 + 3A_1'' +4Q^3 A_1)' +3a_1^2 Q^{(3)}; 
\]
and
\[
G_2 := 6Q^2 (1+A_1)^2 +(6Q^2B_1(1+A_1))' -\frac 12 a_1(9A_1' +3B_1'' +4Q^3 B_1) + \frac 32 a_1^2 Q^4.
\]
Note that both elements, $F_2$ and $G_2$ are in the Schwartz class, but only $F_2$ has the structure of a derivative function of another Schwartz function. Note additionally that $F_2$ and $G_2$ depend on $A_1$, $B_1$ and $a_1$, which are already known.

\medskip

In an independent work (see \cite{MMcomplu}), Martel and Merle discuss the solvability of this linear system. After several computations, they prove the following result.

\begin{lem}
\eqref{Omega2} has a solution $(a_2,A_2,B_2) $ which satisfies $A_2\in S(\R)$, $B_2 =b\varphi +\hat B_2$, with $b<0$ and $\hat B_2 \in S(\R)$. In other words, $B_2\in L^\infty(\R) \backslash L^2(\R)$.
\end{lem}

Note that the new term $A_2(y) Q_{c}^2(y_c)$ is in some sense a \emph{phantom term}, since after the interaction (namely, for $t\sim T_c$), it has almost no size.  Indeed, we obtain
\[
v_4(-T_c,x) \sim Q(y) +Q_c(y-2) +b(Q_c^2)'(y_c),
\]
and
\[
v_4(T_c,x) \sim Q(y) +Q_c(y+2) -b(Q_c^2)'(y_c).
\]
Note that the term $(Q_c^2)'$ is supported near the small soliton.

\medskip

Now we explain the main result in Martel-Merle's paper: the function $B_2$ represents a \emph{defect} in the interaction. Let us explain why. 

\medskip

First of all, note that the term $(Q_c^2)'$ cannot be obtained from the natural shift and scaling variations on the big and small solitons. Indeed, if we have that the shift on the small soliton at time $t=T_c$ is $\Delta_c,$ then
\[
Q_c(y_c +\Delta_c) = Q_c(y_c) + \Delta_c Q_c'(y_c) + \frac12 \Delta_c^2 Q_c''(y_c) +\cdots
\]
However, we cannot obtain the term $(Q_c^2)'$ using the derivatives of $Q_c:$
\[
Q_c'' =cQ_c -Q_c^4,
\]

\[
Q_c^{(3)} =cQ_c' -(Q_c^4)',
\]
and so on...In other words, the term $(Q_c^2)'$ does not appear. Similarly, if $\ve=\ve(y_c)$ is the variation of the scaling of the big soliton, we should have
\[
Q_{1+\ve}(y) = Q(y) + \Lambda Q (y)\ve(y_c) + O_{H^1}(\ve^2),
\]
but every element above is very localized in the $y$ variable (note that $B_2$ is just bounded). The same idea works for the case of the shift on the big soliton.

\medskip

Finally, the best way to understand why $B_2$ represents a defect is by looking at the integrable cases. Note that the previous construction does not depend on the integrability of the equation, so it can be performed for the cases $p=2$ and $p=3$. For instance, if $p=3$, we have
\[
Q_c'' =cQ_c -Q_c^3,
\]

\[
Q_c^{(3)} =cQ_c' -(Q_c^3)',
\]
so we expect that $b(p=3) =0$, which is indeed the case.

\begin{lem}
If $p=3$, one has $B_2 \in \mathcal S(\R)$.
\end{lem}

For the case of the KdV equation ($p=2$), the interpretation is more subtle. We have that 
\[
Q_c^{(3)} =cQ_c' -(Q_c^2)',
\]
so in fact $b\neq 0$, but it represents the third order variation of the shift on the small soliton, in the sense that 
\[
Q_c(y_c +\Delta_c) = Q_c(y_c) + \Delta_c Q_c'(y_c) + \frac12 \Delta_c^2 Q_c''(y_c) + \frac16 \Delta_c^3 Q_c^{(3)}(y_c) +\cdots,
\]
{\bf and} $b=-\frac16 \Delta_c^3$ (note that this coincidence is the important point).

\medskip

\begin{figure}
\begin{center}
\begin{tikzpicture}[
	>=stealth',
	axis/.style={semithick,->},
	coord/.style={dashed, semithick},
	yscale = 1,
	xscale = 1]
	\newcommand{\xmin}{0};
	\newcommand{\xmax}{10};
	\newcommand{\ymin}{0};
	\newcommand{\ymax}{6};
	\newcommand{\ta}{3};
	\newcommand{\fsp}{0.2};
	\draw [axis] (\xmin-\fsp,0) -- (\xmax,0) node [right] {$x$};
	\draw [axis] (0,\ymin-\fsp) -- (0,\ymax) node [below left] {$t$};
	\filldraw[color=light-gray1] (9,1) circle (0.7); 
	\draw [thick,->] (9,1) -- (7,1.5);
	\filldraw[color=light-gray2] (5,1) circle (0.5); 
	\filldraw[color=light-gray2] (5,4) circle (0.5); 
	\filldraw[color=light-gray1] (1.5,4)  circle (0.7); 
	\draw [thick,->] (1.5,4) -- (-0.5, 4.5);
	\draw (5,2.3) node [above] {$y=x-\al(y_c)\sim 0$};
	\draw (1.5,1) node [above] {$t \sim -T_c$};
	\draw (8.5,3.5) node [above] {$t \sim T_c$};
	\draw [domain=0.1:3] plot (\x, {3.7+ sin(350*\x )/20});
	\draw (1.6,3.7) node [below] {$Q_c - b(Q_c^2)'$};
	\draw [domain=7.5:10.5] plot (\x, {0.8+ sin(350*\x )/20});
	\draw (9,0.8) node [below] {$Q_c + b(Q_c^2)'$};
	\draw (9,1.4) node [above] {$-(1-c)<0$};
		\draw (5,0.8) node [above] {$Q(y)$};
\end{tikzpicture}
\end{center}
\caption{Identification of the defect for the symmetric function $v_4$.}\label{7}
\end{figure}
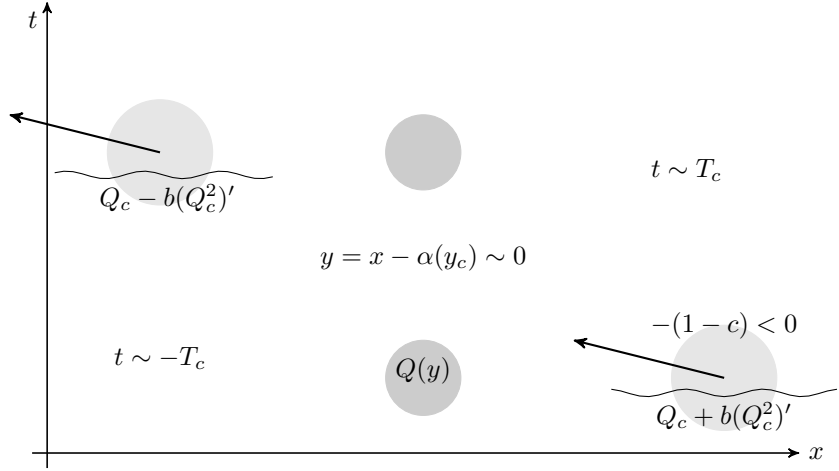

We conclude that every term of the form $b(Q_c^2)'(y_c)$ is \emph{trivial} in the case of the integrable models. However, for the quartic case, it truly represents a \emph{defect} (see Fig \ref{7}).

\medskip

Let us compute the size of $(Q_c^2)'$: from \eqref{solL2} we have
\[
\|(Q_c^2)'\|_{H^1(\R)} \sim c^{2/3 +1/4} =c^{11/12}.
\]
This term will explain later the one in \eqref{3p4}, because it is the first nontrivial element that appears from the interaction. However, in order to prove this fact, we need several additional improvements. First of all, note that after solving \eqref{Omega2}, the largest term in $S[v_4]$ is of the form 
\[
F(y) Q_c^3(y_c).
\] 
Even if we assume that $F \in \mathcal S(\R)$ (the best case), we will have
\[
\|F(y) Q_c^3(y_c)\|_{H^1(\R)} \sim c e^{-\ga_0\sqrt{c} |t|}, \quad \ga_0>0.
\]
Following the stability principle in Claim \ref{SP}, we will have that 
\[
\|v(t) - v_4(t) \|_{H^1(\R)} \lesssim c^{1/2}.
\]
However, $c^{1/2}\gg c^{11/12}$, in other words, we will loose all control on the defect.

\medskip

In \cite{MMcol1} Martel and Merle continue  solving even more complicated linear systems, however this time they use a general theory for solvability. An important problem for that theory is the control of polynomially growing solutions on $y$ (recall that $B_2$ is just bounded, and any equation involving the term $(\mathcal L B)'$ has no right hand side in a derivative form, so it may lead to large solutions).  They use an approximate solution $v_5$ to solve up to order $Q_c^5$, so that now
\[
S[v_5] \sim F(y) Q_c^5,\quad F\in S(\R).
\] 
Moreover, $v_5$ is good enough, in the sense that its higher order terms do not grow too fast in the $y$ variable. Note that 
\[
\|S[v_5](t)\|_{H^1(\R)} \lesssim c^2e^{-\ga_0 \sqrt{c} |t|},
\]
so we will have
\[
\|v(t) -v_5(t)\|_{H^1(\R)} \lesssim c^{3/2} \ll c^{11/12}.
\]
Now the defect becomes evident for the dynamics. 

\medskip

\begin{figure}
\begin{center}
\begin{tikzpicture}[
	>=stealth',
	axis/.style={semithick,->},
	coord/.style={dashed, semithick},
	yscale = 1,
	xscale = 1]
	\newcommand{\xmin}{0};
	\newcommand{\xmax}{10};
	\newcommand{\ymin}{0};
	\newcommand{\ymax}{6};
	\newcommand{\ta}{3};
	\newcommand{\fsp}{0.2};
	\draw [axis] (\xmin-\fsp,0) -- (\xmax,0) node [right] {$x$};
	\draw [axis] (0,\ymin-\fsp) -- (0,\ymax) node [below left] {$t$};
	\filldraw[color=light-gray1] (9,1) circle (0.7); 
	\draw [thick,->] (9,1) -- (7,1.5);
	\filldraw[color=light-gray2] (5,1) circle (0.5); 
	\filldraw[color=light-gray2] (5,4) circle (0.5); 
	\filldraw[color=light-gray1] (1.5,4)  circle (0.7); 
	\draw [thick,->] (1.5,4) -- (-0.5, 4.5);
	\draw (5,2.3) node [above] {$y=x-\al(y_c)\sim 0$};
	\draw (1.5,1) node [above] {$t \sim -T_c$};
	\draw (8.5,3.5) node [above] {$t \sim T_c$};
	\draw [domain=0.1:3] plot (\x, {3.7+ sin(350*\x )/20});
	\draw (1.6,3.7) node [below] {$Q_c - 2b(Q_c^2)'$};
	\draw (9,0.8) node [below] {$Q_c$};
	\draw (9,1.4) node [above] {$-(1-c)<0$};
	\draw (5,0.8) node [above] {$Q(y)$};
\end{tikzpicture}
\end{center}
\caption{Identification of the true defect for the non symmetric function $v_6$.}\label{8}
\end{figure}
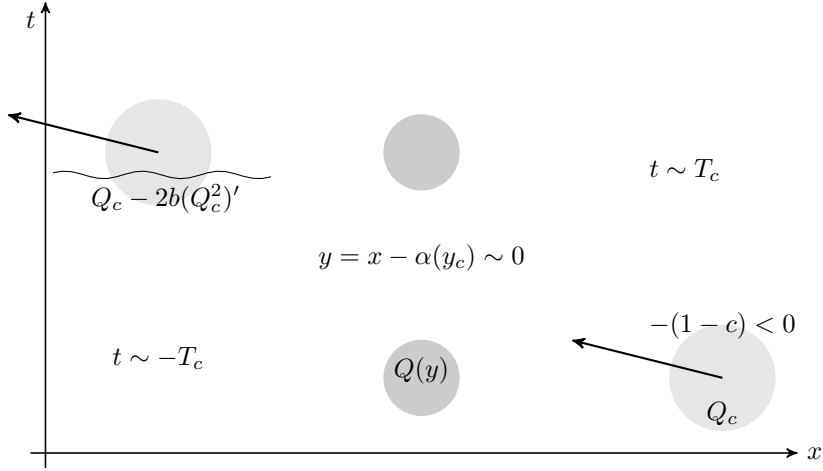

Recall that up to this moment, we have just constructed an approximate solution which has a defect at both sides of time $t=-T_c$ and $t=T_c$. In that sense, this approximate solution is  \emph{symmetric}. We actually need a solution which is almost pure at $t=-T_c$ and not pure at $t=T_c$. In order to obtain such a solution,  one has to modify $v_5$ as follows:
\bee
v_6(t,x) &:=&  Q(y) + Q_c(y_c) +A_1(y) Q_{c}(y_c) +B_1(y) Q_{c}'(y_c) \\
& &   +A_2(y) Q_{c}^2(y_c) +B_2(y) (Q_{c}^2)'(y_c) -b(1+\bar v(y)) (Q_{c}^2)'(y_c)+ \cdots 
\eee
(see Fig. \ref{8}).

Here $\bar v \in \mathcal S$ (the Schwartz class) and it satisfies
\[
B_2(y)  -b(1+\bar v(y))  \sim 0 \quad \hbox{  for $t =-T_c$}, 
\]
and
\[
B_2(y)  -b(1+\bar v(y))  \sim -2b \quad \hbox{  for $t =T_c$}.
\]
Of course we will loose some accuracy on the approximate solution. In that sense the function $\bar v(y)$ allows to loose the minimum degree of accuracy.  In \cite{MMcol1} Martel and Merle proved that 
\[
\|S[v_6](t)\|_{H^1(\R)} \sim c^{3/2} e^{-\ga_0 \sqrt{c} |t|},
\]
so now we have
\be\label{MMF}
\|v(t) -v_6(t)\|_{H^1(\R)} \lesssim c,
\ee
which is still better than $c^{11/12}.$ The defect is now recovered and \eqref{3p4} is proved.  

\medskip

Let us explain more in detail the stability principle announced in Claim \ref{SP}. The idea is to prove the following

\begin{prop}
Assume \eqref{mTe} and $\tilde v$ approximate solution of \eqref{Sv} such that 
\be\label{Res0}
\|\tilde v(-T_c) - Q(\cdot ) -Q_c(\cdot - (1-c)T_c)  \|_{H^1(\R)} \lesssim c^{10},
\ee
and for all time $t\in [-T_c , T_c]$, and $\theta>\frac13$,
\[
\|S[\tilde v](t)\|_{H^1(\R)} \lesssim c^{\theta} T_c^{-1}.
\]
Then the are $\rho(t) \in \R$ such that 
\be\label{Res1}
\|V(t, \cdot -\rho(t)) - \tilde v(t)   \|_{H^1(\R)} \lesssim c^{\theta},
\ee
where $V(t,x):= U(t,x+t)$. Additionally, $\rho'(t)$ is small.
\end{prop}

(In order to prove \eqref{MMF} Martel and Merle take $\theta =1$.) The proof of this result goes as follows. At least, for a certain amount of time $t>-T_c$, one has that \eqref{Res1} is satisfied, mainly because of \eqref{Res0}, \eqref{mTe} and the continuity of the gKdV flow.

\medskip

The idea is to extend this property up to time $t=T_c$. In order to prove this fact, we assume that no matter what is $\rho(t)$, the maximal time for which \eqref{Res1} is satisfied is just $T^*<T_c$. Under this assumption, we will perturb a little bit $V(t)$ by a particular choice of shift $\rho_0(t)$ for which
\be\label{Qp}
\int_\R Q'(y) z(t,x)dx =0, 
\ee
where 
\be\label{zzz}
z(t,x) := V(t, x+ \rho_0(t)) - \tilde v(t,x).
\ee
It is not difficult to see that \eqref{Qp} can be ensured via the Implicit Function Theorem.  Moreover, one has
\[
|\rho_0'(t)| \lesssim \|z(t)\|_{H^1(\R)},
\]
with a constant not depending on $z(t)$. Now we will bootstrap \eqref{Res1}. Indeed, consider the Lyapunov functional
\[
\mathcal F(t) := \frac12\int_\R z_x^2 +\frac12\int_\R (1+\al'(y_c))z^2 -\frac15\int_\R [(\tilde v +z)^5 -\tilde v^5 - 5 \tilde v^4 z].
\]
Note that the term $\al'(y_c)$ is small compared with the constant 1, and it is needed since from \eqref{YAL} one has $y=x-\al(y_c)$, and
\[
\al(y_c) \sim \int_\R Q_c \sim c^{-1/6},
\]
which is a very large perturbation of the soliton center. It is not difficult to see that $\mathcal F$ satisfies the lower bound, uniform on $z$,
\be\label{CoerQQ}
\mathcal F(t) \geq \ga_0 \|z(t)\|_{H^1(\R)}^2 -\frac1{\ga_0}  \abs{\int_\R Q(y) z}^2,
\ee
mainly because $\tilde v(t) \sim Q(y) + \hbox{ small terms}$. On the other hand, we compute the derivative of $\mathcal F$. First, we have that $z(t)$ satisfies the equation
\be\label{ecuZ}
z_t + (z_{xx} -z+ (\tilde v+z)^4 -\tilde v^4)_x +S[\tilde v]  - \rho_0'(\tilde v +z)_x =0.
\ee
Now,
\bee
\mathcal F'(t) & =& \int_\R z_t (-z_{xx} +(1+\al'(y_c))z - (\tilde v +z)^4 + \tilde v^4) \\
& & + \frac12(1-c)\int_\R \al''(y_c) z^2 - \int_\R \tilde v_t [(\tilde v +z)^4 -\tilde v^4 - 4 \tilde v^3 z].
\eee
We replace \eqref{ecuZ} to obtain
\bea
& &\mathcal F'(t) = \int_\R    (z_{xx} -z+ (\tilde v+z)^4 -\tilde v^4)_x    (z_{xx} - (1+\al'(y_c))z + (\tilde v +z)^4 - \tilde v^4)  \nonu\label{11}\\
& & \quad - \int_\R S[\tilde v] (-z_{xx} +(1+\al'(y_c))z - (\tilde v +z)^4 + \tilde v^4) \label{12}\\
& & \quad - \rho_0'  \int_\R (\tilde v +z)_x(z_{xx} -(1+\al'(y_c))z + (\tilde v +z)^4 - \tilde v^4)\nonu \label{13} \\
& & \quad + \frac12(1-c)\int_\R \al''(y_c) z^2 - \int_\R \tilde v_t [(\tilde v +z)^4 -\tilde v^4 - 4 \tilde v^3 z].\nonu \label{14}
\eea
The key estimate above is the one for \eqref{12}. We have
\[
|\eqref{12}| \lesssim \|S[\tilde v](t)\|_{H^1(\R)} \|z(t)\|_{H^1(\R)} \lesssim c^{\theta}T_c^{-1} \|z(t)\|_{H^1(\R)}. 
\]
Note that this estimate is good enough since it depends only linearly on $ \|z(t)\|_{H^1(\R)}$ and not quadratically! For the other terms, the idea is to get estimates of the form
\[
T_c^{-1-\delta_0}\| z(t)\|_{H^1(\R)}^2,
\]
for some $\delta_0>0$, or other better estimates, see \cite{MMcol1} for a detailed proof (at this point the term $\al'(y_c)$ is needed). Under these circumstances, we will have
\[
|\mathcal F'(t)| \lesssim c^{\theta}T_c^{-1} \|z(t)\|_{H^1(\R)} + T_c^{-1-\delta_0}\| z(t)\|_{H^1(\R)}^2,
\]
which, after integration, leads to the bound
\[
|\mathcal F(t)| \lesssim |\mathcal F(-T_c)| + c^{\theta} \sup_{t}\|z(t)\|_{H^1(\R)}  +o_{c\to 0}(1) \sup_t \| z(t)\|_{H^1(\R)}^2.
\]
Using \eqref{CoerQQ}, we will obtain
\[
\sup_{t }\|z(t)\|_{H^1(\R)}^2 \lesssim c^{2\theta}  + \sup_t   \abs{\int_\R Q(y) z}^2,
\]
with constants independent of $z(t)$. Finally, the linear term $\int_\R Qz$ above can be estimated using the conservation of mass, as in the first section of these notes, which leads to \eqref{Res1}. 

\medskip

Finally, some words about \eqref{3p5}. The proof uses an argument by contradiction. If \eqref{3p5} were not true at infinity in time, we can use an stability argument, as the one showed in Section 2, but this time backwards in time, for the sum of two solitons of size 1 and $c$, plus a small error term. The idea is to prove that the size of such an error term is preserved up to the time $t=T_c$. However, note that the constant involved in the coercivity property of Theorem 2 depends on the size of the smaller of both solitons, which in this case is $\sqrt{c}$. This implies that by using the stability argument  backwards near the small soliton we will loose $c^{1/2}$ of accuracy (recall that $c^{11/12+1/2} =c^{17/12}$), leading to a bound of the form (see \eqref{3p4})
\[
 \|Q(y) + Q_c(y_c) - v_6(T_c)\|_{H^1(\R)} \leq \al c^{11/12},
\] 
for any small constant $\al>0$.  However, this bound contradicts the existence of a defect of size $\sim c^{11/12}$.

\bigskip

\section{Stability of particular soliton structures. The case of breathers}

\medskip

In this last chapter we will discuss some very recent results by Miguel A. Alejo and myself. We will place ourselves in an integrable setting, in particular, we will consider the modified Korteweg-de Vries (mKdV) equation 
\be\label{mKdV}
u_t+(u_{xx}+u^3)_x=0, \quad u(t,x)\in \R, \quad (t,x)\in\R^2. 
\ee
Recall that in the case of real-valued initial data, the associated Cauchy problem for \eqref{mKdV} is globally well-posed for initial data in $H^s(\R)$, for any $s> \frac14$, see Kenig-Ponce-Vega \cite{KPV}, and Colliander, Keel, Staffilani, Takaoka and Tao \cite{CKSTT}. Additionally, the (real-valued) flow map is not uniformly continuous if $s<\frac 14$ \cite{KPV2}. In order to prove this last result, Kenig, Ponce and Vega considered a very particular class of solutions of \eqref{mKdV} called \emph{breathers}, discovered by Wadati in \cite{W1}.

\begin{figure}[!htb]
\centering
\includegraphics[width=12.5cm,height=7cm]{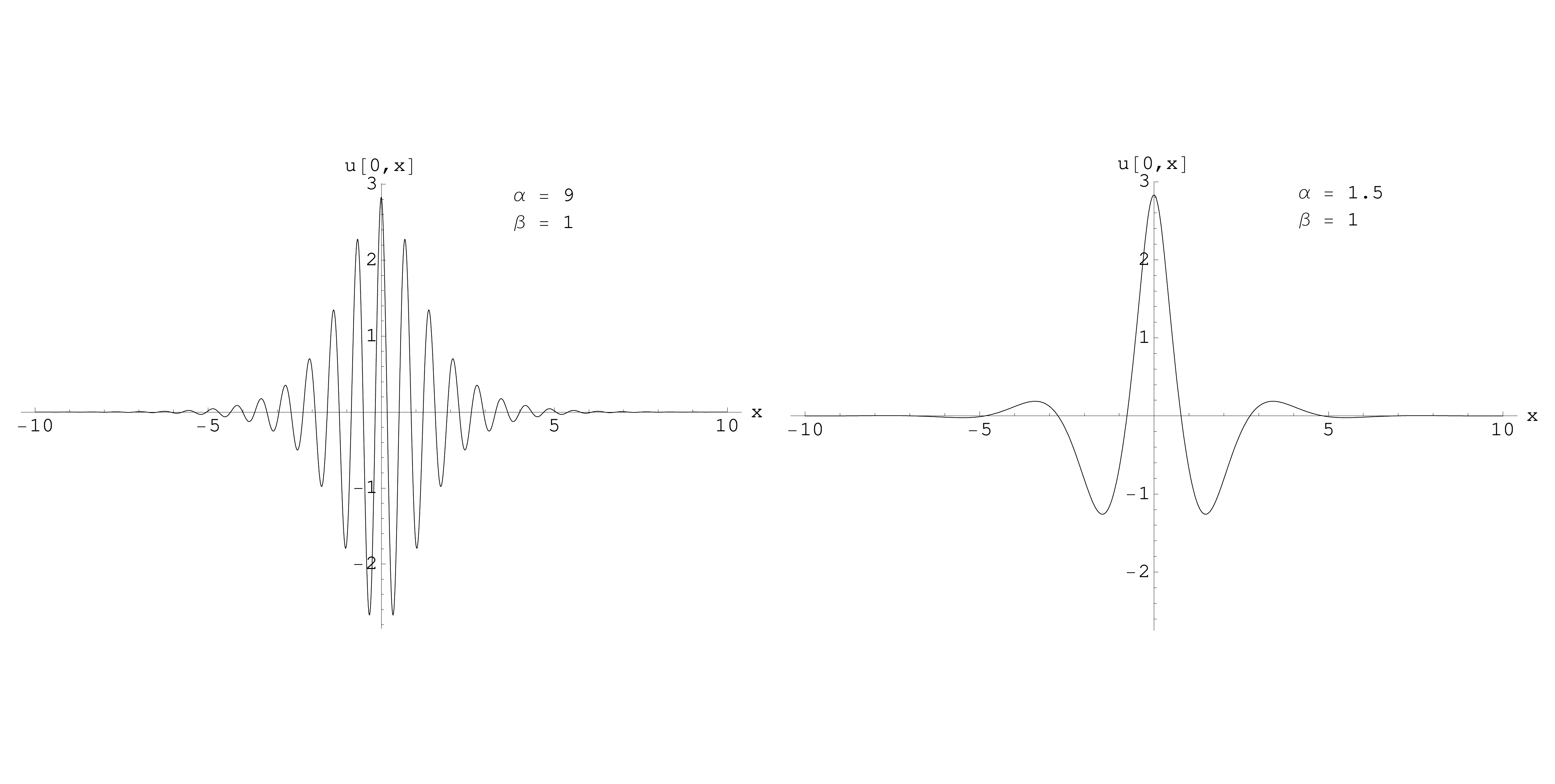}
\small{\caption{Left: mKdV breather \eqref{breather} with $\alpha=9,\beta=1$ at $t=0$.~
Right: mKdV breather \eqref{breather} with $\alpha=1.5,\beta=1$ at $t=0$.\label{BreatherAlfas}}}
\end{figure}

\begin{defn}[See e.g. \cite{W1,La}]\label{mKdVB} 
Let $\al, \bt  >0$ and $x_1,x_2\in \R$ be fixed parameters. The mKdV breather is a smooth solution of \eqref{mKdV} given explicitly  by the formula
\be\label{breather}
\begin{split}
 B:= B(t,x; \al,\bt, x_1,x_2)   := & 2\sqrt{2} \partial_x \Big[ \arctan \Big( \frac{\bt}{\al}\frac{\sin(\al y_1)}{\cosh(\bt y_2)}\Big) \Big] ,
\end{split}
\ee
where
\be\label{y1y2}
y_1 := x+ \delta t + x_1, \quad y_2 : = x+ \ga t + x_2,
\ee
and
\be\label{deltagamma}
\delta := \al^2 -3\bt^2, \quad  \ga := 3\al^2 -\bt^2.
\ee
(See Fig. \ref{BreatherAlfas}.)
\end{defn}

\begin{figure}[!htb]
\centering
\includegraphics[width=11.0cm,height=7.5cm]{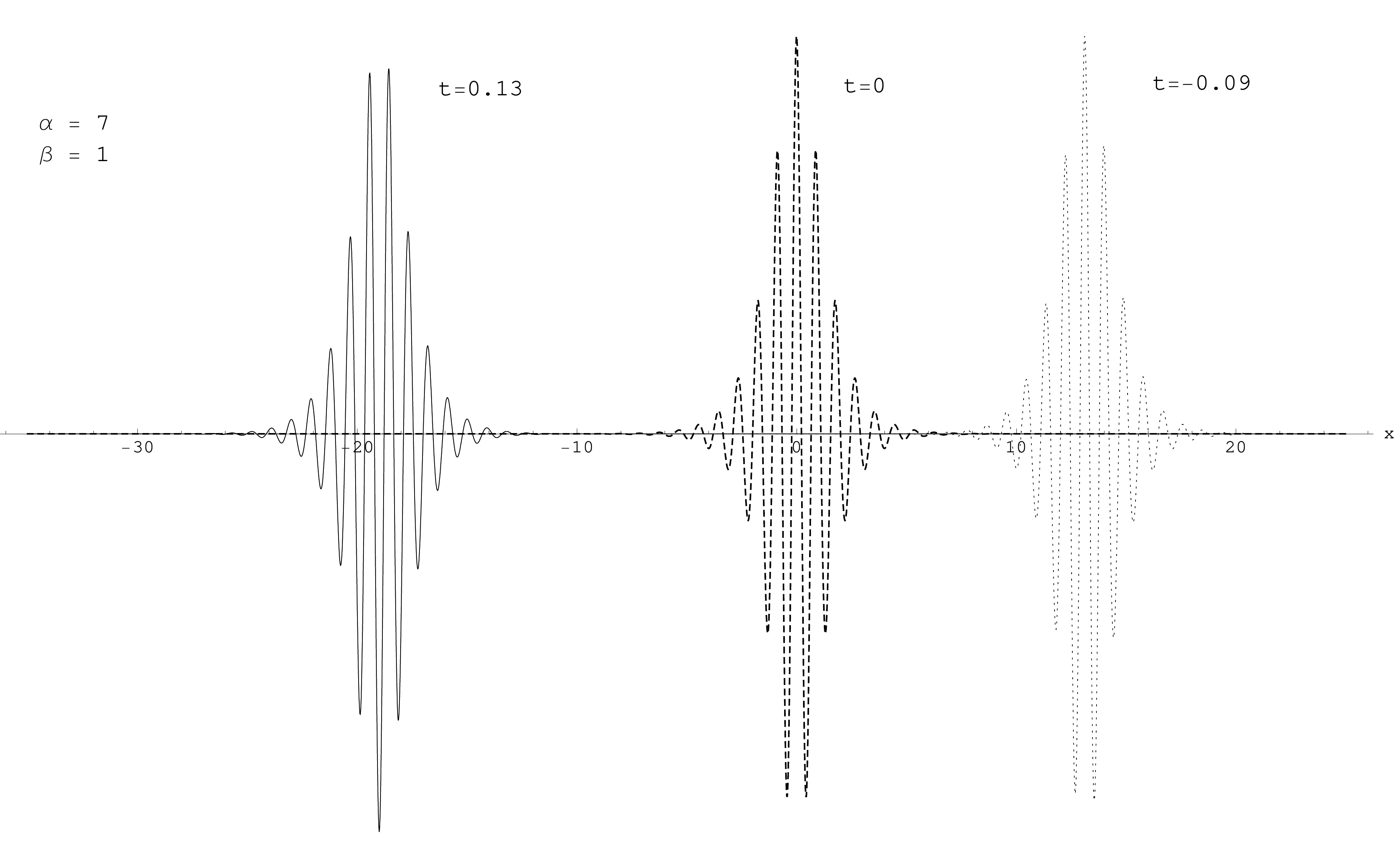}
\small{\caption{ Evolution of the mKdV breather  \eqref{breather} with $\alpha=7,\beta=1$ at instants $t=-0.09$, $t=0$, and $t=0.13$. Note that with the selected values of $\al,\bt$, the \emph{velocity} is given by $\ga=3\al^2 -\bt^2=146>0$ and then the breather moves to the left. (Images taken from \cite{AM}.)\label{EvolBreather}}}
\end{figure}

Breathers are \emph{oscillatory bound states}. They are periodic in time (after a suitable space shift) and localized in space. The parameters $\al$ and $\bt$ are scaling parameters, $x_1,x_2$ are shifts, and $-\ga$ represents the \emph{velocity} of a breather (see Fig. \ref{EvolBreather}).  For a detailed account of the physics of breathers see e.g. \cite{La,AC,Au,Ale,AM} and references therein.

\medskip

Even if the equation is completely integrable in nature, with suitable results on the evolution of well-prepared initial data (see Kruskal et al. \cite{Ga}, Lax \cite{LAX1}, Schuur \cite{Sch}, among others), no rigorous result on the stability of these solutions was given. Numerical computations (see Gorria-Alejo-Vega \cite{AGV}) showed that breathers are \emph{numerically} stable. With Alejo, we give a detailed description of the dynamics around a breather solution. First of all, we proved the following

\begin{thm}[\cite{AM}]\label{Thm4} Breathers are $H^2$-stable. More precisely, for any $\al,\bt>0$, $x_1^0,x_2^0\in \R$, there are $C_0>0$ and $\eta_0>0$ such that for all $\eta\in (0,\eta_0)$ the following holds. Assume that $u_0\in H^2(\R)$ satisfies
\[
\|u_0 - B(0,\cdot ; \al,\bt,x_1^0,x_2^0)\|_{H^2(\R)} <\eta.
\]
then there are $x_1(t)$, $x_2(t) \in \R$ such that 
\[
\sup_{t\in \R} \|u(t) - B(t,\cdot ; \al,\bt,x_1(t),x_2(t))\|_{H^2(\R)} < C_0 \eta.
\]
\end{thm}

The proof of this result is in essence a variational one: we profit of the fact that breathers satisfy very special elliptic equations. 
For similar results in the case of KdV soliton solutions, see e.g. Lax \cite{LAX1} and Maddocks-Sachs \cite{MS}.

\subsection{Sketch of proof of Theorem \ref{Thm4}} Since mKdV is an integrable equation, it has infinitely many conserved quantities. For the proof of Theorem \eqref{Thm4} we will need the additional $H^2$-conserved quantity
\be\label{FFF}
F[u](t) := \frac12\int_\R u_{xx}^2(t) -\frac52\int_\R u^2u_x^2(t) +\frac14\int_\R u^6(t)=F[u](0),
\ee

\begin{figure}
\begin{center}
\begin{tikzpicture}[
	>=stealth',
	axis/.style={semithick,->},
	coord/.style={dashed, semithick},
	yscale = 1,
	xscale = 1]
	\newcommand{\xmin}{0};
	\newcommand{\xmax}{10};
	\newcommand{\ymin}{0};
	\newcommand{\ymax}{6};
	\newcommand{\ta}{3};
	\newcommand{\fsp}{0.2};
	\draw [axis] (\xmin-\fsp,0) -- (\xmax,0) node [right] {$H^2$};
	\draw [axis] (0,\ymin-\fsp) -- (0,\ymax) node [below left] {$F$};
	\draw [dashed] (9,1) -- (1.5, 4);
	\draw [dashed] (5,1) -- (7,4.1);
	\draw [dashed] (5.8,0) -- (5.8,2.4);
	\draw [domain=2:8] plot (\x, {2.7 + (\x-4)*(\x-8)/10});
	\draw (8,2.7) node [above] {$F[u]$};
	\draw (2,4) node [above] {$E[u]=$const};
	\draw (4,1) node [above] {$M[u]=$const};
	\draw (5.8,0) node [below] {$u=B$};	
\end{tikzpicture}
\end{center}
\caption{Schematic representation of Theorem \ref{Thm4}. Each breather has a variational characterization; it is a (local) minimizer of the functional $F$ defined in \eqref{FFF} under suitable constraints on its mass and energy.}\label{10}
\end{figure}
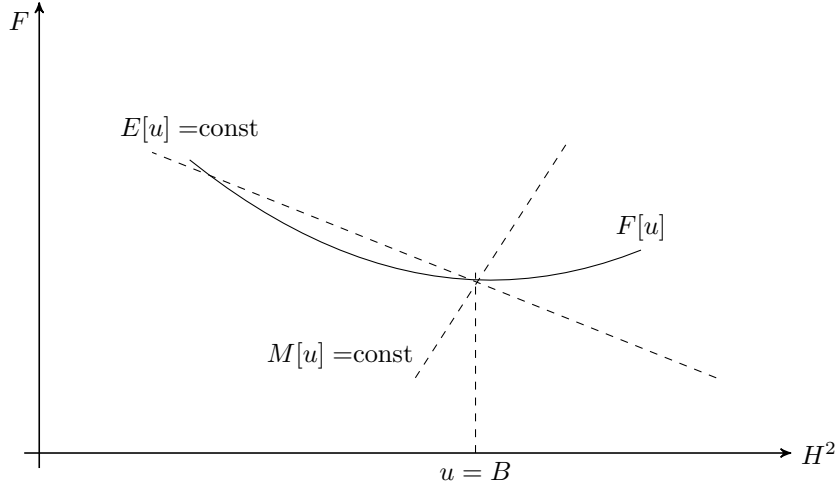

in addition to the standard mass and energy introduced in \eqref{M} and \eqref{E} (take $p=3$). Consider the Lyapunov  functional 
\[
H[u](t) := F[u](t) + 2(\bt^2-\al^2) E[u](t) +(\al^2+\bt^2)^2 M[u](t).
\]
Clearly $H$ is conserved for initial data $u_0\in H^2(\R)$. Moreover, any perturbation of a breather solution of the form 
\[
u(t,x) =B(t,x ; \al,\bt,x_1(t), x_2(t)) +z(t,x),
\] 
with $z$ \emph{small} and $x_1(t)$, $x_2(t)$ to be chosen later, must satisfy the expansion
\be\label{Hu}
H[u](t) = H[B] + \int_\R G[B](t) z + \frac12\int_\R z\mathcal L z + O(\|z(t)\|_{H^2(\R)}^3),
\ee
where $G[B]$ is the nonlinear operator
\bee
G[B] & :=&  B_{4x} -2(\bt^2-\al^2)(B_{xx} +B^3) +(\al^2 +\bt^2)^2 B\\
& &  +5BB_x^2 +5B^2B_{xx} +\frac32B^5,
\eee
and $\mathcal L$ denotes the self-adjoint operator with domain $H^4(\R)$:
\bee
\mathcal L z & := & z_{4x} -2(\bt^2-\al^2)z_{xx} +(\al^2 +\bt^2)^2z  +5B^2 z_{xx} +10B B_x z_x\\
& &  +(5B_x^2 +10BB_{xx} +\frac{15}2 B^4 -6(\bt^2-\al^2)B^2)z.
\eee
One of the key points of the proof is the fact that no matter what are $x_1$ and $x_2$, one has
\be\label{GB}
G[B] \equiv 0.
\ee
In other words, each breather satisfy a suitable fourth order elliptic equation, and in consequence $B$ is a critical point for $H$ (see Fig. \ref{10}). In order to prove this fact, one has two options, either computing \eqref{GB} completely by hand, or proving simpler identities, as is done in \cite{AMV}. First of all, from the definition \eqref{breather} and \eqref{mKdV} we have
\be\label{First}
\tilde B_t +B_{xx} +B^3=0.
\ee
Multiplying this equation by $B_x$ and integrating in space we get
\be\label{2nd}
B_x^2 +\frac12B^4 +2B\tilde B_t -2 \mathcal M_t=0,
\ee
where
\[
\mathcal M := \frac12 \int_{-\infty}^x B^2. 
\]
The third identity that we will need is  the following second order nonlocal equation
\be\label{ide1}
 B_{xt} +  2 \mathcal M_t B   = 2(\bt^2 -\al^2) \tilde B_t  +(\al^2 +\bt^2)^2 B,
\ee
which is the actual equivalent of \eqref{S1} for the case of breathers. This last identity can be proved by hand (see \cite{AM} for a proof). Now we prove \eqref{GB}. We have from (\ref{2nd}) and (\ref{First})
\bee
G[B] & =&  -( B_t + B^3)_{xx} + 2(\bt^2 -\al^2) \tilde B_t  +(\al^2 +\bt^2)^2 B  \\
& & \qquad   + 5 BB_x^2  + 5B^2 B_{xx} + \frac 32 B^5 \\
& =&    - B_{tx} - BB_x^2 + 2B^2B_{xx} + 2(\bt^2 -\al^2) \tilde B_t  +(\al^2 +\bt^2)^2 B + \frac 32 B^5\\
& = & - B_{tx} + B \Big[ \frac 12 B^4 + 2B \tilde B_t -2 \mathcal M_t \Big]  -2B^2 ( \tilde B_t + B^3) + \frac 32 B^5 \\
& &   \qquad + 2(\bt^2 -\al^2) \tilde B_t  +(\al^2 +\bt^2)^2 B\\
& =  & - [ B_{tx} + 2 \mathcal M_t B ] + 2(\bt^2 -\al^2) \tilde B_t  +(\al^2 +\bt^2)^2  B \ = 0.
\eee
In the last line we have used (\ref{ide1}).

\begin{rem}
Some interesting open questions arise from \eqref{GB}. Is $B$ the unique localized solution to \eqref{GB}? If not, under which conditions we recover the uniqueness? Note that for any $x_1,x_2$ shifts, the breather \emph{profile} $B(0,x; \al,\bt,x_1,x_2)$ is solution to \eqref{GB}. In other words, there is a sort of two dimensional set of invariances for \eqref{GB}.
\end{rem}
 
\medskip

From the previous identity it is not difficult to show that 
\[
\mathcal L \partial_{x_1}B =\mathcal L \partial_{x_2}B =0.
\] 
The fact that these two directions span the entire $L^2$ kernel of $\mathcal L$ is not difficult to check. Indeed, the equation $\mathcal L z=0$ is a fourth order ODE whose solutions are spanned by four linearly independent  solutions. Their asymptotic at positive infinity are given by the forms
\[
e^{\pm \bt x} \sin (\al x), \quad e^{\pm \bt x} \cos (\al x),
\] 
so unless $\bt =0$, which is impossible by hypothesis, we only have two linearly independent localized solutions, which coincide with $\partial_{x_1}B $ and $\partial_{x_2}B $ above. 

\medskip

On the other hand, using the Weyl's theorem we have that the continuum spectrum is given by the intervals $[(\al^2+\bt^2)^2,\infty)$ if $\bt\geq \al$, and $[4\al^2\bt^2,\infty)$ if $\bt\leq \al$. 
 
\medskip

We finally consider the problem of counting the number of negative eigenvalues. This is not an easy task, mainly because we deal with a fourth order ODE. The idea is to use the work by L. Greenberg \cite{Gr}, who shows that  the number of negative eigenvalues,
\[
\# \hbox{ negative eigenvalues} = \sum_{x\in \R} \dim \ker W[\partial_{x_1}B ,\partial_{x_2}B ](t,x),
\]
where $W[\partial_{x_1}B ,\partial_{x_2}B ]$ is the Wronskian matrix associated to $\partial_{x_1}B$ and $\partial_{x_2}B$.  The best way for understanding this identity is by considering the same problem for the case of a soliton solution. We have in this case
\[
\# \hbox{ negative eigenvalues} = \sum_{x\in \R} \dim \ker Q' (x),
\]
where we understand $Q'$ as a linear operator in one dimension. Since $Q'(x)=0$ only for $x=0$, we have that  $ \sum_{x\in \R} \dim Q' (x) $ is finite and equals one (i.e., we have just one negative eigenvalue), as expected from Lemma \ref{L1p4}.

\medskip

After some tedious computations we will obtain 
\[
\det W[\partial_{x_1}B ,\partial_{x_2}B ](t,x) =  \frac{16\al^3\bt^3 (\al^2+\bt^2)(\al \sinh(2\bt y_2) -\bt \sin (2\al y_1))}{(\al^2+\bt^2 +\al^2\cosh(2\bt y_2) -\bt^2 \cos(2\al y_1))^2},
\]
so we have nontrivial kernel if and only if 
\be\label{ecy2}
\al \sinh(2\bt y_2) = \bt \sin (2\al y_1).
\ee
Given $x_1$, $x_2$, $t$, $\al$ and $\bt$ fixed, there is only one point $x=x_0\in \R$ (depending on the previous parameters) for which this last identity holds. Indeed, fix $x_1$, $x_2$, $t$, $\al$ and $\bt$. We will look for $\tilde y_2$ solution of 
\be\label{tildey2}
\sinh(\tilde y_2) = \frac{\bt}{\al} \sin \Big(\frac{\al}{\bt} \tilde y_2+  \tilde x_{12} \Big),
\ee
where
\[
\tilde x_{12}:=  2\al [ (\delta-\ga)t + x_1-x_2], \quad \tilde y_2 :=2\bt y_2.
\]
(Recall that $y_1 =x+\delta t+x_1$, $y_2= x+\ga t + x_2$.) If $|\tilde y_2|$ is large enough, say $> M$, there is no solution for this equation.  If now $|\tilde y_2|\leq M$, note that the function
\be\label{fMM}
[-M,M] \ni \tilde y_2 \mapsto \sinh(\tilde y_2) - \frac{\bt}{\al} \sin \Big(\frac{\al}{\bt} \tilde y_2+  \tilde x_{12} \Big)  \in \R
\ee
changes its sign on $[-M,M]$, $M$ large, so it has a root. Moreover, if $\tilde x_{12}\neq 2 k\pi$, $k\in \Z$, such a root is unique (\eqref{fMM} has positive derivative).   Now, if $\tilde x_{12}= 2 k\pi$ for some $k\in \Z$, we will have
\[
 \sinh(\tilde y_2) = \frac{\bt}{\al} \sin \Big(\frac{\al}{\bt} \tilde y_2 \Big),
\]
for which $\tilde y_2 =0$ is a root. If there is another one, it must be unique, by the same reason as before. We conclude that there a unique root $\tilde y_2$ of \eqref{tildey2}. In particular, there is unique $x_0$ satisfying \eqref{ecy2}.

\medskip

At the point $x_0$  we will have
\[
1\leq \dim \ker W[\partial_{x_1}B ,\partial_{x_2}B ](t,x_0) \leq 2,
\]
but it is easy to see that the dimension cannot be 2 since the Wronskian matrix at that point is never identically zero.

\medskip

We conclude that $\mathcal L$ has a unique negative eigenvalue. It is possible to prove that such an eigenvalue is always far from zero, uniformly in time. Even better, following the ideas for the proof of \eqref{coer3}, we are able to prove that there is $c_0>0$ only depending on $\al$ and $\bt$ such that for all $z\in H^2(\R)$, if
\[
\int_\R  z\partial_{x_1}B  =\int_\R z \partial_{x_2}B =0, 
\]
then 
\be\label{5p4}
\int_\R z \mathcal L z \geq c_0 \|z\|_{H^2(\R)}^2 - \frac1{c_0}\abs{\int_\R B z}^2.
\ee
Using the conservation of mass we can estimate the last term above:
\bee
\abs{\int_\R B z(t)} & \leq & \abs{\int_\R B z(0)} + \sup_{t\geq 0} \|z(t)\|_{L^2(\R)}^2 \\
&  \leq&   \|z(0)\|_{L^2(\R)} + \sup_{t\geq 0} \|z(t)\|_{L^2(\R)}^2.
\eee
Replacing in \eqref{5p4}, and using \eqref{Hu}, we conclude.

\begin{rem}
The ideas behind Theorem \ref{Thm4} are very robust and allow to prove different stability (instability) results, provided the linear problem satisfies the desired spectral properties. For example, the sine-Gordon equation (SG)
\be\label{SG}
u_{tt} -u_{xx} + \sin u =0, \quad (u,u_t)(t,x) \in \R^2,
\ee
has a breather solution $(B,B_t)$ (here $B_t$ represents the time derivative of $B$), whose explicit definition is not necessary for these notes.   After some work, we were able to show that $(B,B_t)$ satisfy the elliptic system of equations
\be\label{EcB2}
B_{txx} + \frac 18 B_t^3 +\frac 38 B_x^2 B_t -\frac 14 B_t \cos B -a B_t -\frac b2 B_x =0 ,
\ee
and
\bea\label{EcB}
B_{(4x)} +\frac 38 B_x^2 B_{xx} +\frac 34 B_t B_{tx} B_x +\frac 38 B_t^2 B_{xx} +\frac 58 B_x^2\sin B  -\frac 54 B_{xx}\cos B  & &  \nonu \\
+\frac 14 \sin B\cos B -\frac 18 B_t^2 \sin B   -a(B_{xx} -\sin B) -\frac b2  B_{tx}=0, \qquad 
\eea
for some well-defined constants $a,b\in \R$.  Additionally, there is an associated Lyapunov functional that control the dynamics for all time.  See \cite{AM2} for more details on these ideas. 
\end{rem}

\subsection{$H^1$ stability} It turns out that the previous result can be improved to the level of allowing $H^1$ perturbations. However, now the proof is not variational, since the $H^1$-stability will be a consequence of a \emph{dynamical rigidity} for small perturbations of breather solutions associated to the integrability of the equation.

\begin{thm}\label{Thm5}
Breathers are $H^1$ stable, i.e. stable in the energy space.
\end{thm}

\noindent
{\bf Sketch of proof.} In order to prove this result, we need several preliminary definitions.  First of all, we introduce the complex-valued mKdV soliton.  Consider parameters $\al,\bt>0$, $ x_1$ and $ x_2 \in \R$. Let
\be\label{tQ}
\widetilde Q :=\widetilde Q(x; \al,\bt, x_1, x_2) :=  2\sqrt{2}  \arctan \big( e^{ \bt y_2 + i \al y_1 }\big),
\ee
where $y_1$ and $y_2$ are (re)defined as
\be\label{y1y2mod}
y_1 : = x+ x_1, \quad y_2 := x + x_2.
\ee
We denote the complex-valued soliton profile by
\be
Q :=  \partial_x \widetilde Q  = \frac{2\sqrt{2}  (\bt+i\al) e^{ \bt y_2 + i \al y_1 } }{ 1+e^{ 2(\bt y_2 + i \al y_1 )}} \label{Q0}.
\ee
Finally we denote
\be\label{tQt0} 
\widetilde Q_t :=   -(\bt+i\al)^2 Q.
\ee
We remark that $\widetilde Q$ and $Q$ may blow-up in finite time.  Indeed, assume that 
\be\label{BadCondFF}
\hbox{for $ x_2$ fixed and some $k\in \Z$}, \quad x_{1} =  x_2  +\frac{\pi}{\al} \Big(k+\frac 12 \Big).
\ee
Then $\widetilde Q$ and $Q$ cannot be defined at $x= - x_2.$ However, note that if $x_1 =x_2 =0$, we have $Q(\cdot; \al,\bt,0,0)\in H^1(\R;\Com)$. 

\medskip

Another almost direct conclusion obtained from the definition of $Q$ is the following
\begin{lem}
Fix  $\al,\bt>0$ and $x_1,x_2\in \R$ such that \eqref{BadCondFF} is not satisfied. Then  we have
\be\label{ecQ}
Q_{xx} - (\bt+i\al)^2 Q +Q^3 =0, \quad \hbox{for all } x\in \R,
\ee
and
\be\label{Qx2}
Q_x^2 -(\bt+i\al)^2 Q^2 + \frac 12Q^4 =0, \quad \hbox{for all } x\in \R.
\ee
Moreover, the previous identities can be extended to any $x_1,x_2\in \R$ by continuity.
\end{lem}

Assume that \eqref{BadCondFF} does not hold. Consider the $\sin$ and $\cos$ functions applied to complex numbers. We have from \eqref{tQ} and \eqref{Q0},
\bee
 \sin \Big(\frac{\widetilde Q}{\sqrt{2}}\Big) &=& \sin (2 \arctan e^{\bt y_2 + i\al y_1}) \\
 &=&2e^{\bt y_2 + i\al y_1} \cos^2(\arctan e^{\bt y_2 + i\al y_1})\\
 & = & \frac{ 2e^{\bt y_2 + i\al y_1} }{1+ e^{2(\bt y_2 + i\al y_1)}}=  \frac 1{\bt+i\al}\frac Q{\sqrt{2}}.
\eee
Similarly, from this identity we have
\[
Q_x - (\bt + i\al) \cos \Big( \frac{\widetilde Q}{\sqrt{2}}\Big) Q=0,
\]
so that from \eqref{tQt0} and \eqref{Qx2},
\bee
&  & \widetilde Q_t + (\bt+ i\al)  \Big[  Q_x  \cos \Big(\frac{\widetilde Q}{\sqrt{2}} \Big)   + \frac { Q^2}{\sqrt{2}}\sin \Big( \frac{\widetilde Q}{\sqrt{2}} \Big) \Big]\\
& & \qquad = -(\bt + i\al)^2Q  +    Q_x^2 Q^{-1}   +  \frac 12Q^3  = 0.
\eee
So far, we have proved the following result.

\begin{lem}\label{Equal0}
Let $Q$ be a complex-valued soliton profile with scaling parameters $\al,\bt > 0$ and shifts $x_1,x_2\in \R$, such that \eqref{BadCondFF} is not satisfied. Then we have
\be\label{zero0}
\frac Q{\sqrt{2}}  - (\bt + i\al) \sin \Big(\frac{\widetilde Q}{\sqrt{2}}\Big) \equiv 0 ,
\ee
and
\be\label{zero1}
 \widetilde Q_t    + (\bt+ i\al)  \Big[  Q_x  \cos \Big(\frac{\widetilde Q}{\sqrt{2}} \Big)   + \frac { Q^2}{\sqrt{2}}\sin \Big( \frac{\widetilde Q}{\sqrt{2}} \Big) \Big]  \equiv 0,
\ee
where $\sin z$ and $\cos z$  are defined on the complex plane in the usual sense. 
\end{lem}

\medskip

We introduce now the notion of breather profile. Given parameters  $x_1 ,x_2 \in \R$ and $\al, \bt>0$, we consider $y_1$ and $y_2$ defined in \eqref{y1y2mod}. Let $\tilde B$ be the kink profile
\be\label{tB}
\widetilde B = \widetilde B(x; \al,\bt, x_1, x_2) := 2\sqrt{2} \arctan \Big( \frac{\bt}{\al} \frac{\sin(\al y_1)}{\cosh(\bt y_2)}\Big),
\ee
and with a slight abuse of notation, we redefine
\be\label{BBB}
B:=  \widetilde B_x.
\ee

\medskip

Now we introduce the directions associated to the shifts $ x_1$ and $ x_2$. Given a breather profile of parameters $ \al$, $\bt$, $ x_1$ and $ x_2$, we define
\[
B_1 = B_1 (x; \al, \bt,  x_1, x_2)  := \partial_{x_1} B ,
\]
\[
B_2 =B_2 (x; \al, \bt,  x_1, x_2)  := \partial_{ x_2}  B .
\]
and for $\delta$ and $\gamma$ defined in \eqref{deltagamma}, 
\be\label{tBt}
\widetilde B_t := \delta B_1 +\ga B_2.
\ee
We also have from \eqref{First},
\be\label{tBteq}
\widetilde B_t + B_{xx} + B^3 =0.
\ee

\medskip

We will prove now that there is a deep interplay between complex solitons and  breather profiles. Indeed,

\begin{lem}\label{Equal}
Let $(B,Q)$ be a pair breather-soliton profiles with scaling parameters $\al,\bt>0$ and shifts $x_1,x_2\in \R$. Assume that \eqref{BadCondFF} is not satisfied. Then we have
\be\label{zerot}
\frac{(B-Q)}{\sqrt{2}}  - (\bt-i\al) \sin \Big(\frac{\widetilde B+\widetilde Q}{\sqrt{2}}\Big) \equiv 0 ,
\ee
and
\be\label{zerot2}
\widetilde B_t -\widetilde Q_t    + (\bt-i\al)  \Big[ ( B_x + Q_x) \cos \Big(\frac{\widetilde B+\widetilde Q}{\sqrt{2}} \Big)   + \frac { (B^2 + Q^2)}{\sqrt{2}}\sin \Big( \frac{\widetilde B+\widetilde Q}{\sqrt{2}} \Big) \Big] \equiv 0.
\ee
\end{lem}

\begin{proof}
Let us assume \eqref{zerot} and prove \eqref{zerot2}. We have from \eqref{tQt0} and \eqref{ecQ}
\[
\widetilde Q_t = -(\bt+i\al)^2 Q = -(Q_{xx} + Q^3).
\]
Using \eqref{zerot} we have
\[
B_x-Q_x  - (\bt-i\al)(B+Q) \cos \Big(\frac{\widetilde B+\widetilde Q}{\sqrt{2}}\Big) = 0 ,
\] 
and
\bee
& & B_{xx}-Q_{xx}  - (\bt-i\al) (B_x+Q_x) \cos \Big(\frac{\widetilde B+\widetilde Q}{\sqrt{2}}\Big)\\
& & \qquad + (\bt-i\al) \frac{(B+Q)^2}{\sqrt{2}}\sin \Big(\frac{\widetilde B+\widetilde Q}{\sqrt{2}}\Big) = 0 ,
\eee
so that using once again \eqref{zerot} and \eqref{tBteq} 
\bee
& & \widetilde B_t -\widetilde Q_t    + (\bt-i\al)  \Big[ ( B_x + Q_x) \cos \Big(\frac{\widetilde B+\widetilde Q}{\sqrt{2}} \Big)  
+ \frac { (B^2 + Q^2)}{\sqrt{2}}\sin \Big( \frac{\widetilde B+\widetilde Q}{\sqrt{2}} \Big) \Big] \\
& & = -(B_{xx} + B^3) + Q_{xx} + Q^3   +\Big[  B_{xx}-Q_{xx} + (\bt-i\al) \frac{(B+Q)^2}{\sqrt{2}}\sin \Big(\frac{\widetilde B+\widetilde Q}{\sqrt{2}}\Big) \Big] \\
& &  \qquad +  (\bt-i\al) \frac { (B^2 + Q^2)}{\sqrt{2}}\sin \Big( \frac{\widetilde B+\widetilde Q}{\sqrt{2}} \Big) \\
& & = Q^3-B^3 + \sqrt{2}(\bt-i\al)(B^2 +Q^2 + BQ) \sin \Big(\frac{\widetilde B+\widetilde Q}{\sqrt{2}}\Big) \\
& & =  Q^3-B^3 + (B^2 +Q^2 + BQ)(B-Q)  =0.
\eee 

\smallskip

The proof of \eqref{zerot} is a tedious but straightforward computation which deeply requires the nature of the breather and soliton profiles.  For the proof of this result, see  \cite[Appendix A]{AM1}.
\end{proof}

The previous properties  are consequence of a deeper result. In what follows, we fix a primitive $\tilde f$ of $f$, i.e., 
\be\label{tildeN}
\widetilde f_x := f,
\ee
where $f$ is assumed only in $L^2(\R)$. Notice that even if $f =f(t,x)$ is a solution of mKdV, then a corresponding term $\widetilde f(t,x)$ may be unbounded in space. We introduce the spatial B\"acklund transformation \cite{La}
\be\label{G1}
 G(u_a,u_b,m) := \frac{(u_a-u_b)}{\sqrt{2}}  - m \sin \Big(\frac{\widetilde u_a+\widetilde u_b}{\sqrt{2}}\Big).
\ee
Then
\be\label{BB}
G( Q, 0, \bt+i\al)=0.
\ee
and
\be\label{AAA}
G(B, Q, \bt-i\al) =0.
\ee
There is a second component for this transformation $G$, which involves a time derivative, but for the sake of simplicity we will not deal with it. One can prove (see \cite{AM1}) that given $m$, $u_a$ is solution of mKdV provided $u_b$ is another solution.  For similar examples of the use of B\"acklund transformations, see e.g. the works by Mizumachi and Pelinovski \cite{MP}, and Hoffman and Wayne \cite{EW}.

\medskip

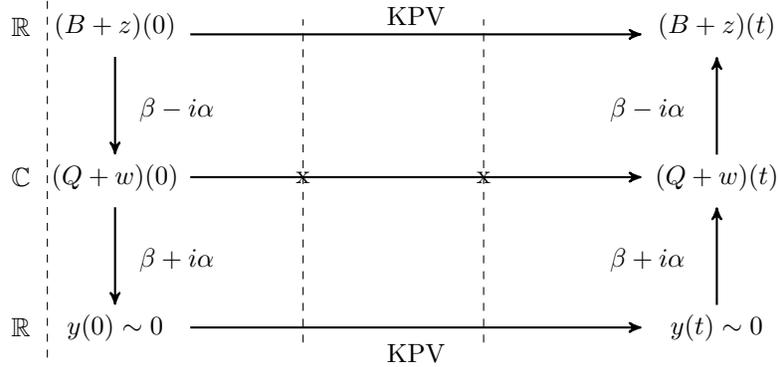
\begin{figure}
\begin{center}
\begin{tikzpicture}[
	>=stealth',
	axis/.style={semithick,->},
	coord/.style={dashed, semithick},
	yscale = 1,
	xscale = 1]
	\newcommand{\xmin}{0};
	\newcommand{\xmax}{10};
	\newcommand{\ymin}{0};
	\newcommand{\ymax}{6};
	\newcommand{\ta}{3};
	\newcommand{\fsp}{0.2};
	\draw (1,4) node [above] {$(B+z)(0)$};
	\draw (1,2) node [above] {$(Q+w)(0)$};
	\draw (1,0) node [above] {$y(0)\sim 0$};
	\draw (9,4) node [above] {$(B+z)(t)$};
	\draw (9,2) node [above] {$(Q+w)(t)$};
	\draw (9,0) node [above] {$y(t)\sim 0$};
	\draw [thick,->] (1,3.9) -- (1,2.6);
	\draw [thick,->] (1,1.9) -- (1,0.6);
	\draw [thick,->] (9,2.6) -- (9,3.9);
	\draw [thick,->] (9,0.6) -- (9,1.9);
	\draw [thick,->] (2,0.3) -- (8,0.3);
	\draw [thick,->] (2,2.3) -- (8,2.3);
	\draw [thick,->] (2,4.2) -- (8,4.2);
	\draw (5,4.2) node [above] {KPV};
	\draw (5,0.2) node [below] {KPV};
	\draw (0,2.3) node [left] {$\Com$};
	\draw (0,4.3) node [left] {$\R$};
	\draw (0,0.3) node [left] {$\R$};
	\draw [dashed] (0.1,-0.1) --(0.1,4.7);
	\draw [dashed] (3.5,0.1) -- (3.5,4.4);
	\draw [dashed] (5.9,0.1) -- (5.9,4.4);
	\draw (3.5,2.3) node {x};	
	\draw (5.9,2.3) node {x};	
	\draw (1.2,3.2) node [right] {$\bt -i\al$};	
	\draw (1.2,1.2) node [right] {$\bt +i\al$};	
	\draw (8.7,3.2) node [left] {$\bt -i\al$};	
	\draw (8.7,1.2) node [left] {$\bt + i\al$};	
\end{tikzpicture}
\end{center}
\caption{An oversimplified description of the proof of the $H^1$ stability of mKdV breathers, using and inverting the B\"acklund transformation twice with parameters $\bt-i\al$ and $\bt+i\al$. See \cite{AM1} for a rigorous proof. KPV means Kenig, Ponce and Vega \cite{KPV}, the symbols $\R$ and $\Com$ mean that the corresponding dynamics (the horizontal lines) are real-valued and complex-valued respectively, and the cross signs are the points where $Q$ blows up (see \eqref{BadCondFF}). Note that the fact that $y(0) \in H^1(\R)$ is real-valued is a nontrivial property of the dynamics, consequence of the fact that $(B+z)(0)$ is real-valued, for $z(0)\in H^1(\R)$ small enough.}\label{11a}
\end{figure}

From the previous paragraph, we see that  the $H^1$ stability result will be a consequence of the Implicit Function Theorem applied to suitable $H^1$ neighborhoods of the points involved in \eqref{BB} and \eqref{AAA}, namely $(B,Q)$ and $(Q,0)$. Note that we are relating the breather to the zero solution via two B\"acklund transformations (see Fig \ref{11a}).

\medskip

The proof of this last result is involved, since $(i)$ we have to prove two invertibility theorems, one near $B$ and another near $Q$, $(ii)$ we deal with initial data which is real valued, but the proof naturally introduces some complex-valued data, however, at  the end we \emph{have} to obtain real-valued data, $(iii)$ this last property is a deep consequence of the particular form of the equation, and $(iv)$ we must be careful with the times where $Q$ blows up. We refer to the reader to \cite{AM1} for a detailed proof.

\bigskip
\bigskip

Some final remarks. The previous results apply without important modifications to the case of the sine-Gordon (SG) equation \eqref{SG} and its corresponding breather \cite[p. 149]{La}. See \cite{BMW,SW,EFM} and references therein for related results.  Since the proofs are very similar, and in order to make this survey non redundant, we skip the details.

\end{document}